\documentclass[12pt]{article}
\usepackage{graphicx}
\usepackage{amsmath,amsfonts,amsthm,amssymb}

\begin{document}
\title
{Small knots and large handle additions}
\author{Ruifeng Qiu and Shicheng Wang
\thanks{Both authors are supported by  NSFC }}
\maketitle \pagestyle{myheadings}

\markboth{R.F. Qiu and S.C. Wang }{Small knots and large handle
additions}

\abstract{We construct a hyperbolic 3-manifold $M$ (with $\partial
M$ totally geodesic) which contains no essential closed surfaces,
but for any even integer $g> 0$ there are infinitely many
separating slopes $r$ on $\partial M$ so that $M[r]$, the
3-manifold obtained by attaching 2-handle to $M$ along $r$,
contains an essential separating closed surface of genus $g$ and
is  still hyperbolic. The result contrasts sharply with those
known finiteness results for the cases $g=0,1$. Our 3-manifold $M$
is the complement of a simple small knot in a handlebody.}
\endabstract

{\bf \S 1. Introduction.} All manifolds in this paper are
orientable. All submanifolds  are embedded and proper ($F\subset
M$ is {\it proper} if $F\cap
\partial M=\partial F$), unless otherwise
specified. A connected  1-manifold (an arc or a circle) on a
surface $F$ is {\it non-trivial} if it does not separate a disc
from $F$.

Let $M$ be a compact 3-manifold with the boundary $\partial M\ne
\emptyset$, $F$ be a surface in $M$ which is not the 2-sphere
$S^2$. Say $F$ is {\it incompressible} if a circle $c\subset F$
bounds a disk in $M$ implies that $c$ bounds a disc in $F$. Say a
surface in $M$ is {\it essential} if either it is incompressible
and is not parallel to a sub-surface  of $\partial M$, or it is a
2-sphere which does not bound a 3-ball in $M$. Say a 3-manifold
$M$ is {\it irreducible} if each 2-sphere in $M$ bounds a 3-ball.
Say $M$ is {\it $\partial$-irreducible} if $\partial M$ is
incompressible. Say $M$ is atoroidal if it contains no essential
tori; Say $M$ is {\it anannular} if it contains no essential
annuli.

Say a 3-manifold $M$ is {\it simple} if $M$ is irreducible,
$\partial$-irreducible, anannular and atoroidal. Suppose $M$ is a
simple 3-manifold with $\partial M\ne \emptyset$. By Thurston's
theorem,  $M$  admits a complete finite volume hyperbolic
structure with totally geodesic boundary (with torus components in
$\partial M$ removed) [T]. A knot $K$ in $M$ is {\it simple} if
$M_K$, the complement of $K$ in $M$, is simple. A 3-manifold $M$
is {\it small} if $M$ contains no essential closed surface. A knot
$K$ in $M$ is {\it small} if $M_{K}$ is small.

A {\it (separating) slope} $r$ in $\partial M$ is the isotopy
class of a non-oriented non-trivial (separating) circle in
$\partial M$.  We denote by $M[r]$ the manifold obtained by adding
a 2-handle to $M$ along a regular neighborhood of $r$ in $\partial
M$ and then capping off spherical components with 3-balls.
Specially, if $r$ lies in a torus component of $\partial M$, this
operation is known as Dehn filling.

Essential surface is a basic tool to study 3-manifolds and handle
addition is a basic method to construct 3-manifolds. A central
topic connecting  those two aspects in 3-manifold topology is the
following:

\vskip 0.3 true cm

 {\bf Question 1.} {\it Let $M$ be a simple
3-manifold with $\partial M\ne \emptyset$ which contains no
essential closed surface of genus $g$. How many slopes $r\subset
\partial M$ are there so that $M[r]$ contains an essential closed surface
of genus $g$?}

{\bf Remark on Question 1.} The mapping class group of a simple
3-manifold $M$ with $\partial M\ne \emptyset$ is finite. The
question is asked only for simple 3-manifolds  to avoid possibly
infinitely many slopes produced from  Dehn twists along essential
discs or annuli.

The main result in this paper is the following:

\vskip 0.3 true cm

 {\bf Theorem 1.} {\it There is a  simple small
knot $K$ in the handlebody $H$ of genus 3 such that for any even
integer $g> 0$, there are infinitely many separating slopes $r$ in
$\partial H$ so that $H_{K}[r]$ contains an essential separating
closed surface of genus $g$. Moreover those $H_K[r]$ are still
simple.}

{\bf Remarks on  Theorem 1.}

(1) Suppose $M$ is a simple 3-manifolds with $\partial M\ne
\emptyset$.

(i) $\partial M$ is a torus. W. Thurston's pioneer result claims
that at most finitely many slopes on $\partial M$ so that $M[r]$
are not hyperbolic [T], hence the number of slopes are finite in
 Question 1 when $g=0, 1$. The sharp upper bound of  such
slopes are given by C. Gordon and J. Luecke and by Gordon when
$g=0,1$, see [G] for a survey. A. Hatcher proved the number of
slopes in Question 1 is finitely many  for all $g$ [H].

(ii) $\partial M$ has genus $>1$.  M. Scharlemann and Y-Q Wu [SW]
have shown that if $g=0, 1$, then there are only finitely many
separating slopes $r$ so that $M[r]$ contains an essential closed
surface of genus $g$. Very recently, M. Lackenby [L] generalized
Thurston's finiteness result to handlebody attaching, that is to
add 2-handles simultaneously. He proved that for a hyperbolic
3-manifold $M$ there is a finite set $\cal C$ of exceptional
circles on $\partial M$ so that
 attaching a handlebody to $M$
is still hyperbolike if none of those circles is attached to a
meridian disc of the handlebody.

Theorem 1 and those finiteness results of [T], [H], [SW] and [L]
give a globe view about the answer of Question 1. In particular
those finiteness results of [T], [H] and [SW] do not hold in
general. We think the example in Theorem 1 also indicates that the
finiteness result of [L] does not hold in general  (a working
project of the authors).

(2) It is unusual to the authors that a given manifold $M$
provides "non-finiteness" answer to Question 1 for all even genus
$g\ge 2$. From aesthetic point of view, one may wonder if there is
a manifold provides "non-finiteness" answer to Question 1 for all
genus $g\ge 2$. We think that the answer is positive. In this case
the knot $K$ is complicated and then the proof of that $H_K$ is
small will be much more difficult (a working project of the
authors).

(3) Without handle addition, the 3-manifold $M$ itself in Theorem
1 is interesting independently. First the construction of the
small knot in Theorem 1 can be modified to provided infinitely
many small knots in handlebodies of any genus $g>1$ (a working
project of the authors). Up to our knowledge, no examples of
simple small knots in the handlebody of genus $>1$ were explicitly
presented before. Secondly $M$ provides a hyperbolic 3-manifold
with totally geodesic boundary which splits over essential
surfaces of genus $g$ in infinitely many different ways for each
even $g>0$.

{\bf Remarks on the Proof of Theorem 1 and the organization of the
paper.} In \S 2 we construct a knot $K$ and infinitely many
separating surfaces $S_{g,l}$ of genus $g$  for each even $g > 0$
in the handlebody $H$ of genus 3,   such that all those surfaces
$S_{g,l}$ are disjoint from  $K$ and have connected boundaries.
 Those $\partial S_{g,l}$ will be served as the slopes $r$ in
Theorem 1. Some elementary properties of $S_{g,l}$ and of $K$ are
also described in \S 2. Let $\hat S_{g,l}\subset H_K[\partial
S_{g,l}]$ be the closed surface obtained by capping off $\partial
S_{g,l}$ with a disk. In \S 3 we will prove that $\hat S_{g,l}$ is
incompressible in $H_K[\partial S_{g,l}]$ as well as  that
$\partial S_{g, l}$ and $\partial S_{g,l'}$ are not isotopic in
$\partial H$ when $l\ne l'$. \S 4 and \S 5 are devoted to prove
that the knot $K$ is simple and small.

A result in [J] is quoted in \S 3, which is a crucial step for the
proof of Proposition 3.1, and a result in [CGLS] is quoted in \S
4, which is used to shorten the argument of Case 2 in the proof of
Lemma 4.4. Up to those two results and the knowledge in the
beginning of standard textbooks of elementary algebraic topology,
combinatorial groups and 3-manifolds, the paper is self-contained.
Even so, the argument of Case 1 (2) in the proof of Lemma 4.4 is
initialed  by Gordon-Litherland in  middle 1980's.

\vskip 0.5 true cm

{\bf \S 2. Construction of the surfaces $S_{g,l}$ and the knot $K$
in $H$.}

Suppose $X_1$ and $X_2$ are connected proper sub-manifolds  of $M$
with complementary  dimensions and meet tranversely. Let $||X_1,
X_2||$ be the absolute value of their algebraic intersection
number. Since all manifolds are orientable, $||X_1, X_2||$ is well
defined. For a compact manifold $X$, $|X|$ denotes the number of
components of $X$.  If $X$ is an arc or an annulus, we often use
$\partial_1 X$ to denote one component of $\partial X$ and
$\partial_2 X$ to denote another.

Let $H$ be the handlebody of genus 3. Let  $\bigl\{B_{1}, B_{2},
B_{3}\bigr\}$ be a set of basis disks of $H$, and
$\bigl\{E_{1},E_{2}\bigr\}$ be two separating disks of $H$ which
separate $H$ into three solid tori $J_{1}, J_{2}$ and $J_{3}$. See
Figure 2.1.

\vskip 0.3 true cm

\begin{center}%
\includegraphics[totalheight=4cm]{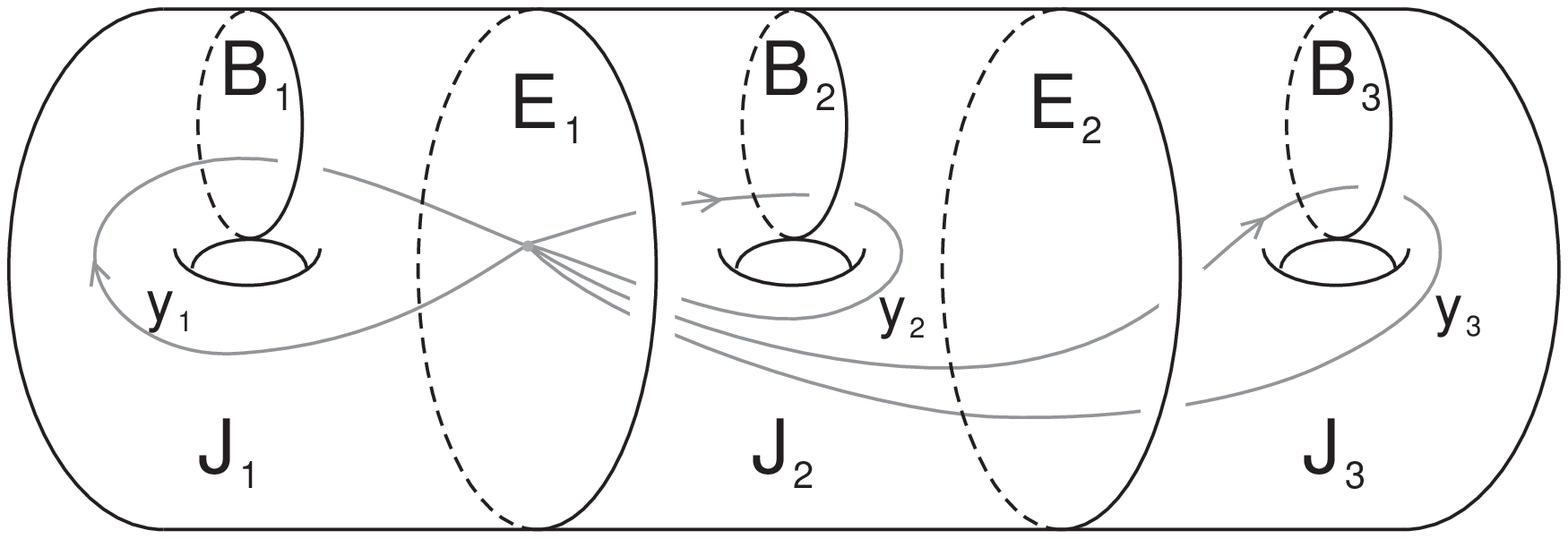}%
\begin{center}%
Figure 2.1
\end{center}
\end{center}

The  orientable surface $S_g$ of even genus $g>0$ with $|\partial
S_g|=1$ can be presented as in Figure 2.2 (where $g=4$). Each
surface $S_{g,l}$ we are going to construct in $H$ can be viewed
as a properly embedded image of $S_g$, where the disk in Figure
2.2 is sent to $E_1$ (approximately) and the 1-handle ended at
$v_i$ and $u_i$ is sent to the 1-handle $N(\alpha_i)$ attached to
$E_1$, which will be shown in Figures 2.3 and 2.4.

\vskip 0.5 true cm

\begin{center}%
\includegraphics[totalheight=4cm]{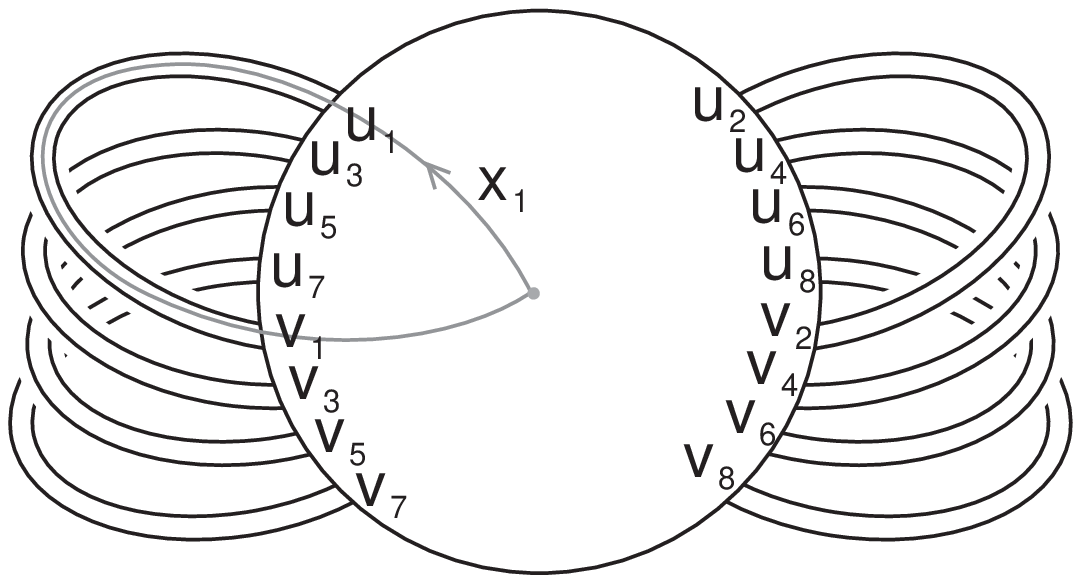}%
\begin{center}%
Figure 2.2
\end{center}
\end{center}

{\bf Remark on Figure 2.2.} In Figure 2.2 if we attach $g$
1-handles on each side of the disc for odd $g$ in the same way, we
get a surface of genus $g-1$ with three boundary components rather
then a surface of genus $g$ with one boundary component.

Let $C$ be a closed curve in $\partial H$ (with one
self-intersection) as in Figure 2.3. Then $\partial E_{1}\cup
\partial E_{2}$ separates $C$ into eight embedded arcs
$c_{1},\ldots,c_{8}$, where $c_{3},c_{7}\subset J_{1}$ with
$||\partial B_{1}, c_7||=3$, $||\partial B_{1}, c_3||=1$;  $c_{2},
c_{4}, c_6, c_8\subset J_{2}$ with $||\partial B_{2}, c_4||=1$,
$||\partial B_{2}, c_6||=3$, $||\partial B_{2}, c_2||=||\partial
B_{2}, c_8||=0$; $c_{1}, c_{5}\subset J_{3}$ with $||\partial
B_{3}, c_1||=3$, $||\partial B_{3}, c_5||=1$.

Let $u_{1},\ldots,u_{2g}, v_{1},\ldots,v_{2g}$ be $4g$ points
located on $\partial E_{1}$ in the cyclic order $u_1$, $u_3$, ...
, $u_{2g-3}$, $u_{2g-1}$,  $v_1$, $v_3$, ... , $v_{2g-3}$,
$v_{2g-1}$, $v_{2g}$, $v_{2g-2}$, ... ,$v_4$,$v_2$, $u_{2g}$,
$u_{2g-2}$, ... ,$u_4$, $u_2$ as in Figure 2.3 (see also Figure
2.2).

\vskip 0.3 true cm
\begin{center}%
\includegraphics[totalheight=4cm]{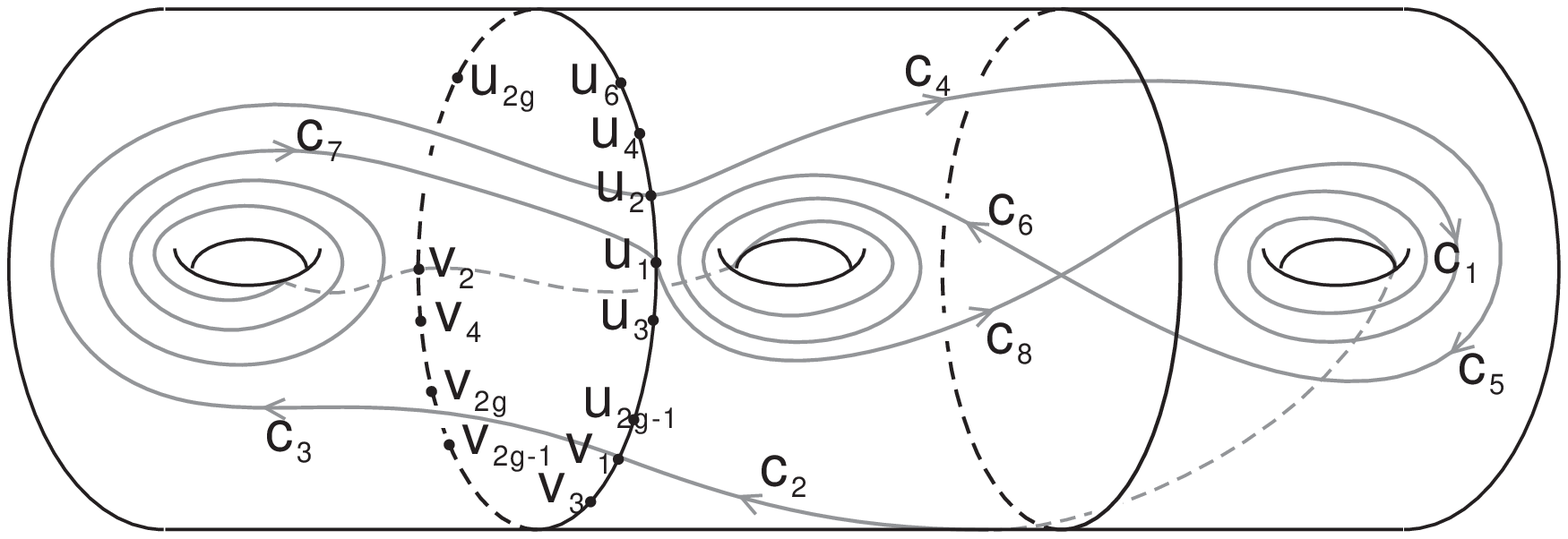}%
\begin{center}%
Figure 2.3
\end{center}
\end{center}

By the order of those points, we can assume that the isotopy has
been made so that $\partial (c_8\cup c_1\cup c_2) =\{u_1, v_1\}$,
$\partial c_3 =\{v_1, u_2\}$, $\partial c_7 =\{v_2, u_1\}$. Then
pick a proper arc $c_*$, (resp. $c_\#$) in $\partial H\cap
(J_2\cup J_3)$ connecting $v_3$ and $v_2$ (resp. $u_2$ and $u_3$)
as in Figure 2.4, where $||c_*,\partial B_2||=l,$ $ l\ge 3$.

Now we  define  oriented arcs on $\partial H$ to connect some
pairs in $\{u_i, v_j;\, i,j=1,...,2g\}$ as follow: First let
$\overline{u_{1}v_{1}}=c_8\cup c_{1}\cup c_2$,
$\overline{v_{1}u_{2}}=c_{3}$, $\overline{v_{2}u_{1}}=c_{7}$,
$\overline{u_{2}u_{3}}=c_\#$, $\overline{v_{3}v_{2}}=c_{*}$. Then
let $\overline{v_{2i}u_{2i-1}}$ and $\overline{v_{2i-1}u_{2i}}$ be
a proper arcs on $\partial H\cap J_1$ parallel to
 $c_{7}$ and  $c_3$ respectively, $i=2,...,g$
 and $\overline{u_{2i}u_{2i+1}}$ and
$\overline{v_{2i+1}v_{2i}}$ be a  proper arcs on  $\partial H\cap
(J_2\cup J_3)$ parallel to $c_\#$ and $c_*$ respectively,
$i=2,...,g-1$. See $\overline {u_3v_4}$ and $\overline {u_4u_5}$
in Figure 2.4. Now define

$\alpha_1$=$\overline{u_1v_1}$,\hskip 8.2truecm (2.0)

and for $1< 4k+j\le 2g, j=1,2,3,4$,

$\alpha_{4k+1}=\overline{u_{4k+1}u_{4k}}\cup\alpha_{4k}\cup
\overline{v_{4k}v_{4k+1}}$,  \hskip 3.3truecm (2.1)

$\alpha_{4k+2}=\overline{v_{4k+2}u_{4k+1}}\cup\alpha_{4k+1}\cup
\overline{v_{4k+1}u_{4k+2}}$,  \hskip 2truecm(2.2)

$\alpha_{4k+3}=\overline{v_{4k+3}v_{4k+2}}\cup\alpha_{4k+2}\cup
\overline{u_{4k+2}u_{4k+3}}$, \hskip 2truecm(2.3)

$\alpha_{4k+4}=\overline{u_{4k+4}v_{4k+3}}\cup\alpha_{4k+3}\cup
\overline{u_{4k+3}v_{4k+4}}$. \hskip 2truecm(2.4)

 Hence $\alpha_{k-1}\subset \alpha_k$ is an increasing sequence
of embedded arcs on $\partial H$.

\vskip 0.3 true cm

\begin{center}
\includegraphics[totalheight=4cm]{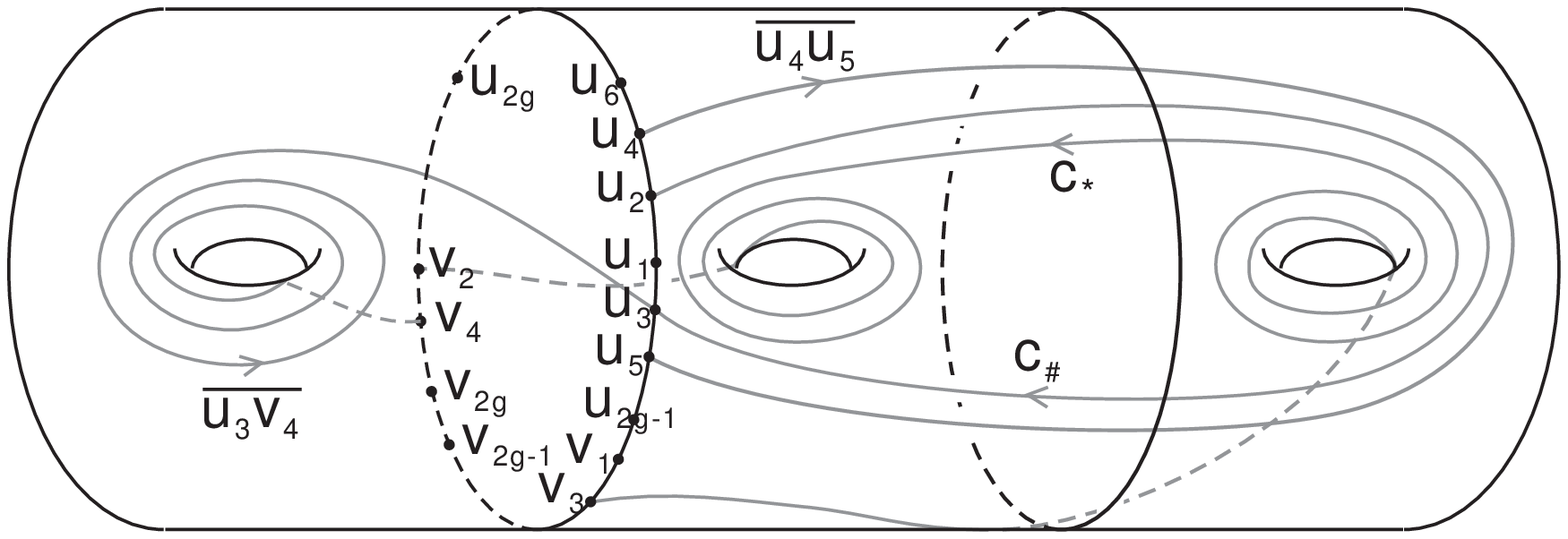}%
\begin{center}
Figure 2.4
\end{center}
\end{center}

Let $\alpha\subset
\partial H$ be an arc which meets $\partial S$ exactly in its two
ends for a proper separating surfaces $S\subset H$. The resulting
proper surface by tubing $S$ along $\alpha$ in $H$, denoted by
$S(\alpha)$, is obtained by first attaching 2-dimensional 1-handle
$N(\alpha)\subset
\partial H$ to $S$, then making  the surface $S\cup N(\alpha)$ to
be proper, that is, pushing the  interior of $S\cup N(\alpha)$
into the interior of $H$. The image of $N(\alpha)$ after the
pushing is still denoted by $N(\alpha)$. Since $S$ is orientable
and separating, it  is a direct observation that $S(\alpha)$ is
still   orientable  and separating.

Since $\alpha_1$ meets $E_1$ exactly in its two ends, we  do
tubing of $E_1$ along $\alpha_1$ to get $E_1(\alpha_1)$. Now
$\alpha_2$ meets $E_1(\alpha_1)$ exactly in its two ends, we do
tubing of $E_1(\alpha_1)$ along $\alpha_2$ to get $E_1(\alpha_1,
\alpha_2)=E_1(\alpha_1) (\alpha_2)$, where the tube $N(\alpha_2)$
is thinner and closer to $\partial H$ so that it goes over the
tube $N(\alpha_1)$. Hence $E_1(\alpha_1, \alpha_2)$ is proper
embedded surface. Repeating this process by tubing along
$\alpha_3,..., \alpha_{2g}$ in order we get a  surface
$E_1(\alpha_1,...., \alpha_{2g})$, denoted by $S_{g,l}$, in $H$.
Clearly $S_{g,l}$ is a proper embedding of $S_g$ into $H$ for each
even $g>0$. We survey this fact as

\vskip 0.3 true cm

{\bf Lemma 2.1.} {\it $S_{g,l}$ is an orientable separating
surface in $H$. Moreover $S_{g,l}$ is  of genus $g$ with
$|\partial S_{g,l}|=1$  for even $g>0$ (and of genus $g-1$ with
$|\partial S_{g,l}|=3$ for odd $g$).}

In the construction of $S_{g,l}$  for all $g,l$, we may assume
that (i) the positions of the arcs $\alpha_1, \alpha_2$ are fixed;
(ii) each tube $N(\alpha_i)$ has distance $\delta/i$ from
$\alpha_i$ for some $\delta>0$. By (i) and (ii), we have (iii)
$N(\alpha_1)$, $N(\alpha_2)$ and the part of $N(\alpha_3)$ goes
over $N(\alpha_2)$ are fixed for all $g,l$.

Now our knot $K$ is obtained by pushing $C$ into the interior of
$H$ along the inward normal direction of $\partial H$ in the
following way:  (iv) first push the arc $c_7\cup c_8\cup c_{1}\cup
c_2\cup c_3$ to stay between $N(\alpha_2)$ and $N(\alpha_3)$, (v)
then push the arc $c_4\cup c_5\cup c_6$ so that it has distance
larger than $\delta/3$ from $\partial H$ and is disjoint from
$N(\alpha_1)$. Below we use $a_i$ to denote the image of $c_i$
after pushing. Then $E_{1}\cup E_{2}$ separates $K$ into 8 arcs
$a_{1},\ldots,a_{8}$.  See Figure 2.5 for $K, a_i\subset H$, where
$a_6$ is crossing under $a_8$, and $||a_i, B_k|| \text{(in
$H$)}=||c_i,
\partial B_k|| \text{(in $\partial H$)}$, $i=1,...,8$ and $k=1,2,3$.
\begin{center}%
\includegraphics[totalheight=5cm]{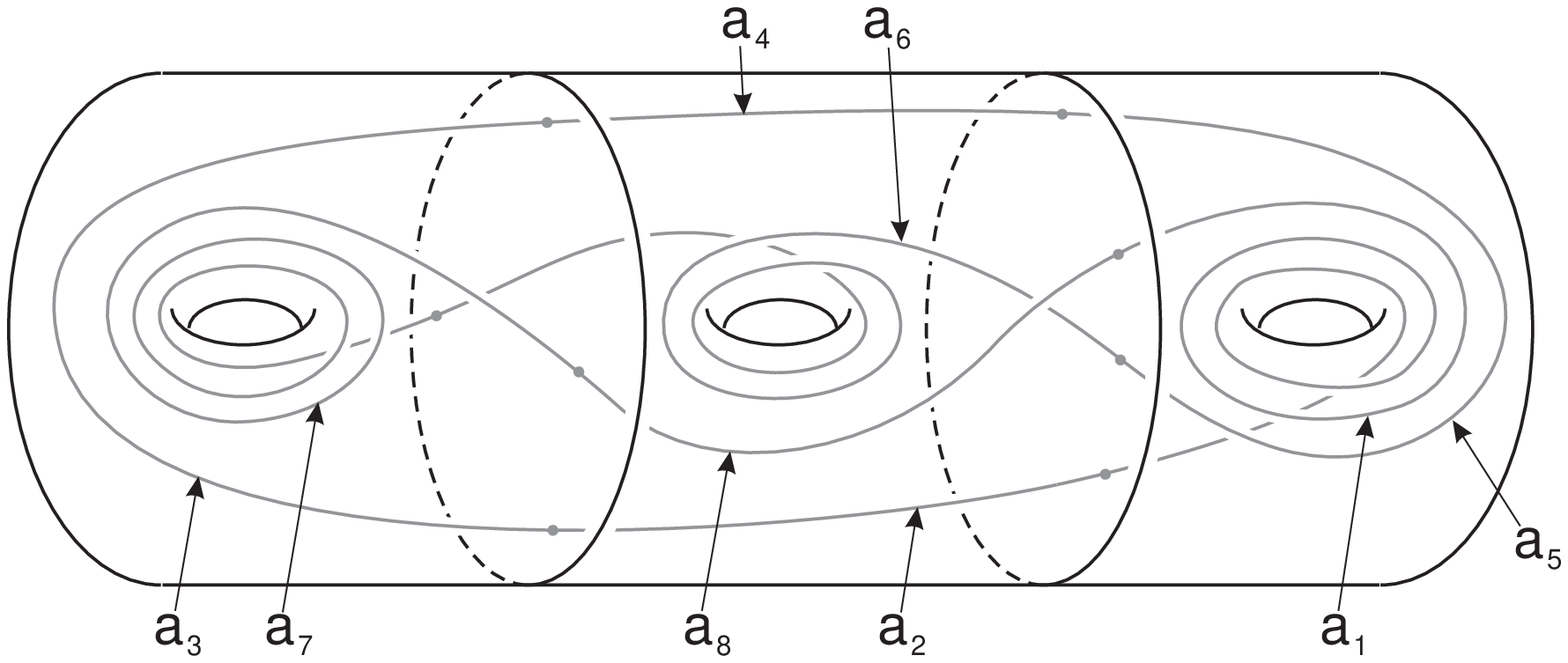}%
\begin{center}%
Figure 2.5
\end{center}
\end{center}

{\bf Lemma 2.2.}  {\it $K\cap S_{g,l}=\emptyset$ for all $g,l$.}

{\bf Proof.}  By (iii) and (iv), the part $a_7\cup a_8\cup
a_{1}\cup a_2\cup a_3$ of $K$ is disjoint from $S_{g,l}$. By (ii),
(iii) and (v), the part $a_4\cup a_5\cup a_6\subset J_2\cup J_3$
of $K$ is also disjoint from $S_{g,l}$. Hence $K\cap
S_{g,l}=\emptyset$ for all $g,l$.\qed

Let $N(K)=K\times D$ be the regular neighborhood of $K$ in $H$
such that (i) $S_{g,l}\subset H_K=H-intN(K)$ for all $g,l$, (ii)
the product structure has been adjusted so that
$\cup_{i=1}^{8}\partial a_i\times D\subset E_1\cup E_2$. Let
$F_{j}=E_{j}-intN(K)$, $j=1,2$; $M_{k}=H_{K}\cap J_{k}$,
$k=1,2,3$; and $T=\partial (K\times D)$. Then $F_{1}\cup F_{2}$
separates $T$ into eight annuli $A_{1},\ldots, A_{8}$, where
$A_{i}=a_{i}\times \partial D$. Moreover $K$ and $C$ bound a
non-embedded annulus in $H$ (the trace of pushing $C$ to $K$)
which is cut by $E_{1}\cup E_{2}$ into eight disk $D_{i*}$, with
$a_i\subset
\partial D_{i*}$, $i=1,...,8$. Suppose $a_i\subset M_k$, then
$D_i=D_{i*}\cap M_k$ is a proper disc in $M_k$. Let
$W_{i}=\overline{\partial N(D_i\cup A_i)-\partial M_k} $, where
$N(D_i\cup A_i)$ is a regular neighborhood of $D_i\cup A_i$ in
$M_k$. Then $W_{i}$ is a proper separating disk in $M_{k}$. Each
$W_{i}$ intersects $F_{1}\cup F_{2}$ in two arcs $l_{i}$ and
$l_{i+1}$. Note $W=\cup_{i\neq 6}W_{i}$ is still a (non-proper)
disc. We use $\mu$ to denote the meridian slope on $T$.  See
Figure 2.6 for $A_i, D_i, W_i, l_i, l_{i+1}\subset M_k$.

\vskip 0.3 true cm

\begin{center}%
\includegraphics[totalheight=4.5cm]{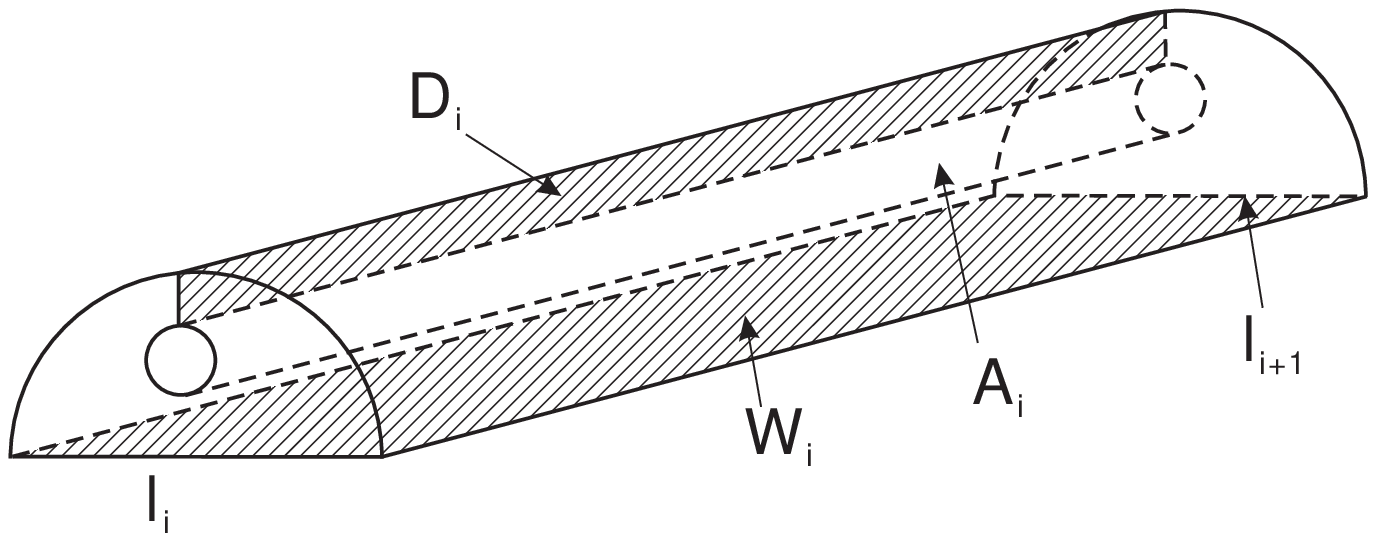}%
\begin{center}%
Figure 2.6
\end{center}
\end{center}

 The following facts about $K$ and $a_i$, which are based on Figure
 2.5 and whose proofs involve only elementary algebraic topology rather
than 3-manifold topology, will be used in \S 4 and \S 5.

\vskip 0.3 true cm

 {\bf Lemma 2.3.} {\it (1) $K$ is not
contractible in $H$.

(2) Suppose  $a_i, a_j\subset J_k$.  There is no relative homotopy
on $ (J_k, E_1\cup E_2)$ which either sends  $a_i$ to $E_1\cup
E_2$; or sends $a_i$ to $a_j$ unless $(i,j)$ is $(2,8)$.

(3) Suppose $a_i, a_j\subset J_k$. The meridians of $A_i$ and
$A_j$ are not homotopic in $M_k$.

(4) Suppose $B$ is a proper disc of $J_k$ with $|B\cap E_j|\le 1$,
$j=1,2$. If $|B\cap (\cup_{a_i\subset J_k} a_i)|< 3-|B\cap
(E_1\cup E_2)|$, then $B$ separates a 3-ball from $J_k$.

(5) There is no annulus $A\subset H$ such  that (i) $\partial_1
A=K$ and $\partial_2 A\subset \partial H$, (ii)  each component of
$A\cap (E_1\cup E_2)$ is non-trivial in $A$.}

{\bf Proof.}  The proofs of (1) (2) (3) are directly.

(4) If $B$ is a separating disk in $J_k$, then $B$ separates a
3-ball from $J_k$, since $J_k$ is a solid torus. So we need only
to show that each non-separating disk in $J_k$ does not meet the
inequality in (4).

Note $B\cap E_j$ is either an arc or empty-set. May suppose $B$ is
non-separating disk in $J_2$ which meets each $E_j$ in an arc
$d_j$, $j=1,2$ (the remaining cases are more directly). Let
$b_{j}$ be an arc in $E_{j}$ connecting the two endpoints of
$a_{6}$ and $a_8$. Then $c=b_1\cup a_6\cup b_2 \cup a_8$ is  an
circle which goes around $J_3$ three times. Hence

$3=||B, c||\le ||B, a_6\cup a_8||+||\partial B, b_1||+||\partial
B, b_2||$

$=||B, a_6\cup a_8||+||d_1, b_1||+||d_2, b_2||.$

By Jordan Curve Theorem, $||d_j, b_j||\le 1=|B\cap E_j|$. Hence

$|B\cap (\cup_{a_i\subset J_k} a_i)|\ge ||B, b_1\cup a_6\cup b_2
\cup a_8||\ge 3-|B\cap (E_1\cup E_2)|$.

(5) Otherwise there is an annulus $A$ meets (i) and (ii) in (5).
Then by (1),  $A$ is cut by $ E_{1}\cup E_{2}$ into eight
rectangles  $R_i$, $i=1,...,8$, each $R_i$ has two opposite sides
in $E_1\cup E_2$ and remaining two sides  $a_i$ in $K$ and
$a_{i}^{*}\subset
\partial H$. Let $b_{i}$  be an arc in $E_{1}\cup E_2$ connecting
$\partial a_{i}^{*}$, and denote the circle $b_{i}\cup a_{i}^{*}
\subset J_k$ by $e_i$, $i=1,3,7,5$. In a basis of $H_1(\partial
J_3,Z)$, $e_{1}$ and $e_{5}$ have coordinates $(3,p)$ and $(1, q)$
respectively, and hence $||b_1, b_5||=||e_1, e_5||= |3q-p|\ne 0$,
since $p$ and $3$ are co-prime. It follows that  $\partial b_1=
\partial a_1^*$ and $\partial b_5=\partial a_5^*$ are alternating
on $\partial E_2$. By the same reason, $\partial a_3^*$ and
$\partial a_7^*$ are alternating on $\partial E_1$.

Now back to $J_2$,  $\partial a_i^*$'s have the cyclic order
$4,8,2,6$ in  $\partial E_1$, and the cyclic order $4, 8, 6, 2$ in
$\partial E_2$. Hence there are four disjoint arcs on $E_1\cup
E_2$ such that the two with $a_4^*\cup a_8^*$ form a circle
$e_{4,8}$ on $\partial J_2$, and the other two with $a_2^*\cup
a_6^*$ form a circle $e_{2,6}$ on $\partial J_2$, moreover
$e_{2,6}$ and $e_{4,8}$ are disjoint, therefore they are parallel
on $\partial J_2$. But in a basis of $H_1(\partial J_2,Z)$, those
two circles have coordinates $(3,p)$ and $(1, q)$, and $||e_{2,6},
e_{4,8}||=|3q-p|\ne 0$, since $3$ and $p$ are coprime. A
contradiction. \qed

 \vskip 0.5 true cm

 {\bf \S 3. Proof of Theorem 1 by  assuming that $K$ is simple and small.}

In this section $g>0$ will be even integer. By Lemma 2.1, let
$\hat S_{g,l}\subset H_K[\partial S_{g,l}]\subset H(\partial
S_{g,l})$ the surface obtained by capping off the boundary of
$S_{g,l}$ with a disk. Then $\hat S_{g,l}$ is a closed surface of
genus $g$.

Now Theorem 1 follows from the following two propositions (the
"Moreover" part of Theorem 1 follows directly from [SW]).

\vskip 0.3 true cm

{\bf Proposition 3.0} {\it $K\subset H$ is a simple and small
knot.}

\vskip 0.3 true cm

{\bf Proposition 3.1.}   {\it (1) $\hat S_{g,l}$ is incompressible
in $H_K[\partial S_{g,l}]$.

(2) for given $g$, $\partial S_{g,l}$ and $\partial S_{g, l'}$ are
not the same slope in $\partial H_K$ when $l\ne l'$.}

We choose the center of $E_1$ as the common base point for the
fundamental groups of $H$ and of all surfaces $S_{g,l}$. Now
$\pi_{1}(H)$ is the free
 group of rank three generated by $y_{1}, y_{2}, y_{3}$ indicated in Figure 2.1.
and  $\pi_{1}(S_{g,l})$ is the free group of rank $2g$ generated
by $x_{1},\ldots,x_{2g}$, where $x_i$ is the generator given by
 $\alpha_i$ and  two arcs in $E_1$ (see Figure 2.2).   Let
$\phi :S_{g,l}\to H$ be the inclusion (precisely $\phi$ should be
$\phi_{g,l}$, we omit sub-index without making confusion),
$\phi_*: \pi_1(S_{g,l})\to \pi_1(H)$ be the induced homomorphism.
It is easy to see from Figures 2.1, 2,3 and 2.4,

$\tilde c_3=y_1,\, \tilde c_7=y_1^3,\,  \tilde c_\#=y_2y_3,\,
\tilde c_*=y_3^{-3}y_2^{-l}$\hskip 5 true cm (*)

where $\tilde c_3\in \pi_1(H)$ is given by  $c_3$ and two arcs in
$E_1$ and so on.

Recall that $\overline{v_{2i}u_{2i-1}}$,
$\overline{v_{2i-1}u_{2i}}$ $\overline{u_{2i}u_{2i+1}}$ and
$\overline{v_{2i+1}v_{2i}}$ are
 parallel copies of
 $c_{7}$, $c_3$, $c_\#$ and $c_*$ respectively. One can read $\phi_*(x_i)$ directly
as words in $y_1, y_2, y_3$ by (2.0)---(2.4) and (*). They are:

$\phi_*(x_1)=y_3^3$, \hskip 10.5 truecm(3.0)

and for $1< 4i+j\le 2g, j=1,2,3,4$,

$\phi_*(x_{4i+1})=y_{3}^{-1}y_{2}^{-1}\phi_*(x_{4i})y_{2}^{l}y_{3}^{3}
=(w_1^{-1}w_2)^iy_3^3 (w_1w_2^{-1})^{i},$ \hskip2.4truecm(3.1)

$\phi_*(x_{4i+2})=y_{1}^3\phi_*(x_{4i+1})y_{1}=y_{1}^3(w_1^{-1}w_2)^iy_3^3
(w_1w_2^{-1})^{i}y_{1},$ \hskip 2.5 truecm(3.2)

$\phi_*(x_{4i+3})=y_{3}^{-3}y_{2}^{-l}\phi_*(x_{4i+2})y_{2}y_{3}
=w_2(w_1^{-1}w_2)^iy_3^3 (w_1w_2^{-1})^{i}w_1$, \hskip
1truecm(3.3)

$\phi_*(x_{4i+4})=y_{1}^{-1}\phi_*(x_{4i+3})y_{1}^{-3}=y_{1}^{-1}w_2(w_1^{-1}w_2)^iy_3^3
(w_1w_2^{-1})^{i}w_1 y_{1}^{-3}$, \hskip 0.5 truecm(3.4)

where  $w_1=y_1y_2y_3$ and $w_2=y_3^{-3}y_2^{-l}y_1^{3}$. \hskip
5.5 truecm(3.5)

Obviousely

$ (w_1w_2^{-1})^{j}w_1(w_1^{-1}w_2)^i =w_2(w_1^{-1}w_2)^{i-j-1}$
if $i>j$ and

$(w_1w_2^{-1})^{j}w_1(w_1^{-1}w_2)^i= (w_1w_2^{-1})^{j-i}w_1$ if
$i\le j$. \hskip 4 truecm (3.6)

$(w_1w_2^{-1})^{j}w_2(w_1^{-1}w_2)^i =w_2(w_1^{-1}w_2)^{i-j}$ if
$i\ge j$ and

$(w_1w_2^{-1})^{j}w_2(w_1^{-1}w_2)^i=(w_1w_2^{-1})^{j-i-1}w_1$ if
$i< j$.  \hskip 3.5 truecm (3.7)

\vskip 0.3 true cm

{\bf Lemma 3.2} {\it (1) $S_{g,l}$ is incompressible in $H$.

(2) for given $g$, $\partial S_{g,l}$ and $\partial S_{g, l'}$ are
not in the same slope in $\partial H$ if $l\ne l'$.}

{\bf Proof}  By (3.5), the right sides of (3.0)---(3.4) are
reduced words in $<y_1, y_2, y_3>$. Now we present
$\pi_1(S_{l,g})$ as the free product $G_1*G_2$, where
$G_1=<x_1,x_3,..., x_{2g-1}>$ and $G_2=<x_2,x_4,..., x_{2g}>$.

(1) We need only to show that $\phi_*: \pi_1(S_{g,l})\to \pi_1(H)$
is injective.

For each $w_2\in G_2$, we may suppose that $w_2$ is a reduced from
in $<x_2,..., x_{2g}>$. Now we can present $\phi_*(w_2)$ as a word
in $<y_1, y_2, y_3>$ by first  replacing each $x_{2i}^{\pm l}\in
w_2$ by $\phi_*(x_{2i})^{\pm l}$ and then applying (3.2) and
(3.4). By (3.5) (3.6), (3.7) and obvious cancellations, one can
get a reduced form of $\phi_*(w_2)$ in $<y_1, y_2, y_3>$. Indeed
by an induction on the length of the reduced form $w_2$, it is
easy to see that if $w_2\ne 1$, then $\phi_*(w_2)\ne 1 $ and
$\phi_*(w_2)$ has the reduced form started from and ended by the
non-zero powers of $y_1$. Similarly  one can argue that for $1\ne
w_1\in G_1$, $\phi_*(w_1)\ne 1 $ and $\phi_*(w_1)$ has the reduced
form started from and ended by the non-zero powers of $y_3$ and
$y_2$.

Now present each $1\ne w\in G_1*G_2$ in a reduced form
$g_1g_2...g_n$ of $G_1*G_2$, and each $g_i$ in a reduced form in
$G_1$ or $G_2$. It is clear that $\phi_*(w)\ne 1$.

(2) For given $g,l$,  the conjugacy class corresponding to
$\partial S_{g,l}$ in $\pi_1(S_{g,l})$ can be  presented by a
reduced word below (see Figure 2.2):

$x_1x_3^{-1}...x_{2g-3}x_{2g-1}^{-1}x_1^{-1}x_3...x_{2g-3}^{-1}x_{2g-1}
x_{2g}^{-1}x_{2g-2}...x_{4}^{-1}x_2x_{2g}x_{2g-2}^{-1}...x_{4}x_2^{-1}$
\hskip 0.5 truecm(**) .

Now we can present $\phi_*([\partial S_{g,l}])$ in $\pi_1(H)$ as a
word of $<y_1, y_2, y_3>$  by (**) and (3.0)---(3.4). Then doing
cancellations to get the reduced form of $\phi_*([\partial
S_{g,l}])$ is very directly and
 all powers
of $y_2$ are untouched in this process. It follows that $
\phi_*([\partial S_{g,l}])$ and $ \phi_*([\partial S_{g,l'}])$ do
not have the same cyclic reduced form when $l\ne l'$. Hence if
$l\ne l'$, $S_{g,l}$ and $S_{g,l'}$ are not homotopic in $H$, and
therefore there are not isotopic in $\partial H$.
 \qed

Now $S_{g,l}$  separates $H$ into two components $P_{1}$ and
$P_{2}$ with $\partial P_{1}=T_{1}\cup S_{g,l}$ and $\partial
P_{2}=T_{2}\cup S_{g,l}$, where $T_1\cup T_2=\partial H$ and
$\partial T_1=\partial T_2=\partial S_{g,l}$.

\vskip 0.3 true cm

{\bf Lemma 3.3.} {\it $T_{i}$ is incompressible in $H$.}

{\bf Proof.} Let $\phi_\# : H_1(S_{g,l},Z)\to H_1(H,Z)$ be the
induced homomorphism  on the first homology groups. Note that
$H_{1}(H, Z)=Z+Z+Z$ is generated are $\bar y_{1}$, $\bar y_{2}$
and $\bar y_{3}$, where $\bar y_i=\pi(y_i)$, and $\pi:
\pi_1(H,Z)\to H_1(H,Z)$ is the abelization. By (3.0)--(3.4), it is
easy to see that $i_{\#} (H_{1}(S_{g,l},Z))$ is a subgroup of
$H_{1}(H,Z)$ generated by $4\bar y_1$, $(l+1)\bar y_{2}$, $\bar
y_3$. Thus $H_{1}(H,Z)/\phi_{\#}(H_{1}(S_{g,l},Z))$ is a finite
group (of order $4l+4$).

If $T_{i}$, $i=1$ or 2, is compressible, then there is a
compressing disk $B'_1$ in $H$ for $T_i$. Since $\partial
B'_1\cap\partial S_{g,l}=\emptyset$ and $S_{g,l}$ is
incompressible in $H$, by standard argument in 3-manifold
topology, we may assume that $B'_1\cap S_{g,l}=\emptyset$.
Furthermore, since $H$ is a handlebody, we may also assume that
$B'_1$ is non-separating in $H$. Thus there are two properly
embedded disks $B'_2$ and $B'_3$ in $H$ such that $(B'_1,B'_2,
B'_3)$ is a set of basis disks of $H$. Let $z_{1}$, $z_{2}$ and
$z_{3}$ be generators of $\pi_{1}(H)$ corresponding to $B'_1,
B'_2$ and $B'_3$. Since $S_{g,l}$ misses $B'_1$,
$\phi_*(\pi_{1}(S_{g,l}))\subset G\subset \pi_{1}(H)$, where $G$
is generated by $z_{2}$ and $z_{3}$. Then clearly
$H_{1}(H,Z)/\phi_\#(H_{1}(S_{g,l},Z))$ is infinite group, a
contradiction.\qed

\vskip 0.3 true cm

 {\bf Jaco's Lemma [J].} {\it Let $M$ be a
compact 3-manifold with compressible $\partial M$ and $r$ be a
circle in $\partial M$. If $\partial M- r$ is incompressible in
$M$, then either $M[r]$ is a 3-ball or $\partial M[r]$ is
incompressible.}

{\bf Proof of Proposition 3.1.} Since $S_{g,l}$ is incompressible
in $H$ by Lemma 3.2 and $H$ contains no closed incompressible
surface, $\partial P_{i}$ is compressible in $P_{i}, i=1,2;$

Since $T_{1}, T_{2}$ and $S_{g,l}$ are incompressible in $H$ by
Lemma 3.3 and Lemma 3.2, $T_{i}$ and $S_{g,l}$ are incompressible
in $P_{i}$. Hence $\partial P_i-\partial S_{g,l}$ is
incompressible in $P_i$, $i=1,2.$ Since  clearly $P_{i}[\partial
S_{g,l}]$ is not a 3-ball,  $\partial P_{i}[\partial S_{g,l}]$ is
incompressible by Jaco's Lemma. It follows that $\hat S_{g,l}$,
which is parallel to a component of
 $\partial P_{i}[\partial S_{g,l}]$, is
incompressible in $P_{i}[\partial S_{g,l}]$. Since $H[\partial
S_{g,l}]$ is a union of $P_{1}[\partial S_{g,l}]$ and
$P_{2}[\partial S_{g,l}]$ along $\hat S_{g,l}$, $\hat S_{g,l}$ is
incompressible in $H[\partial S_{g,l}]$. Therefore $\hat S_{g,l}$
is incompressible in $H_K[\partial S_{g,l}]$. We proved
Proposition 3.1 (1).

 Proposition 3.1 (2) follows Lemma 3.2. (2)\qed

\vskip 0.5 true cm

 {\bf \S 4. $H_k$ is irreducible,
$\partial$-irreducible, anannular.}

Recall $E_j$, $F_j$, $J_k$, $M_k$, $B_k$, $a_i$, $A_i$, $D_i$,
$T$, $\mu$ defined in \S 2.

\vskip 0.3 true cm

{\bf Lemma 4.1.} {\it $F_{1}\cup F_2$ is incompressible and
$\partial$-incompressible in $H_{K}$.}

{\it Proof.} \  Suppose first $F_{1}\cup F_{2}$  is compressible
in $H_{K}$. Then there is a disk $B\subset M_k$ such that $B\cap
(F_{1}\cup F_{2})=\partial B$ and $\partial B$ is a non-trivial
circle on $F_{1}\cup F_{2}$. Denote by $B^{'}$ the disk bounded by
$\partial B$ in $E_{1}\cup E_2$. Then $B\cup B^{'}$ is a 2-sphere
$S^2$ in the solid torus $J_k$,  so $B\cup B^{'}$ bounds a 3-ball
$B^3$ in $J_k$. Since $\partial B$ is non-trivial in $F_{1}\cup
F_2$, $B^{'}$ contains $\partial_1 a_i$ for some $a_i\subset J_k$.
Since $S^2$ is separating and $a_i$ is connected, we must have
$(a_{i},
\partial a_i)\subset (B^3, B')$, which provides a relative
homotopy on $(J_k, E_1\cup E_2)$ sending $a_i$ to $E_1\cup E_2$,
see Figure 4.1 (a), which contradicts Lemma 2.3 (2).

Suppose then $F_{1}\cup F_{2}$ is $\partial$-compressible in
$H_K$. Then there is a non-trivial arc $a$ in $F_{1}\cup F_{2}$
and an arc $b$ in $\partial H_{K}$ bound a proper disk $B$ in
$M_k$. There are two cases:

(1) \ $b\subset T$. Then $b$ is a proper arc in $A_{i}$,
$i=1,5,3,7$. Now either $b$ is a trivial arc in $A_i$, then there
 is an arc $b'$ in $\partial A_i$ such that  the circle $a\cup b'$ is
non-trivial in $F_1\cup F_2$ but bounds a disc in $M_k$, which
contradicts the incompressibility of $F_1\cup F_2$ we just proved;
or $b$ is a non-trivial arc in $A_i$,  then the disc $B$ provides
a relative homotopy on $(J_k, E_1\cup E_2)$ sending $a_i$ to
$E_1\cup E_2$, which contradicts Lemma 2.3 (2).

(2) \ $b\subset \partial H$. Since $|B\cap (E_1\cup E_2)|=1$, $B$
separates a 3-ball $B^3$ from $J_k$  by Lemma 2.3 (4). Since $a$
is non-trivial in $F_{1}\cup F_2$, by the same reason as the end
of the first paragraph, one of $a_{i}$ lies in $B^3$ with
$\partial a_i$ lies in a disc in $\partial B^3\cap E_1$, see
Figure 4.1 (b),  which contradicts Lemma 2.3 (2). \qed

\vskip 0.3 true cm

\begin{center}%
\includegraphics[totalheight=10cm]{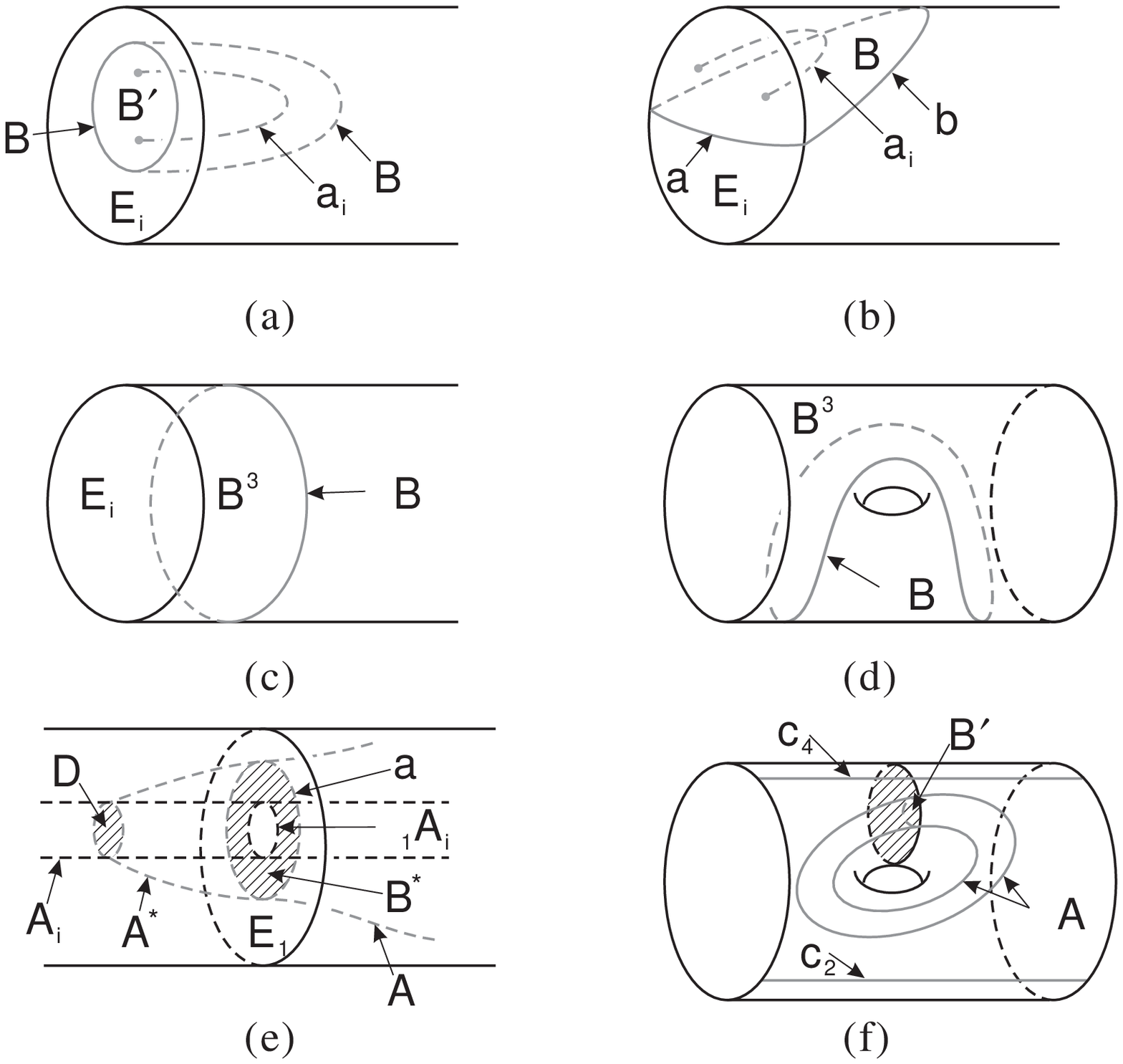}%
\begin{center}%
Figure 4.1
\end{center}
\end{center}

\vskip 0.3 true cm

{\bf Lemma  4.2.} {\it $H_{K}$ is irreducible.}

{\bf Proof.} Otherwise there is an essential 2-sphere $S^2$ in
$H_{K}$. Since $H$ is irreducible, $S^2$ bounds a 3-ball $B^3$ in
$H$ with $K\subset B^3$, which contradicts Lemma 2.3 (1).\qed

\vskip 0.3 true cm

{\bf Lemma 4.3.} {\it $H_{K}$ is $\partial$-irreducible.}

{\bf Proof.} \ Suppose $H_{K}$ is $\partial$-reducible. Let $B$ be
a compressing disk of $\partial H_{K}$. If $\partial B\subset T$,
then $H_K$ contains an essential 2-sphere, which contradicts Lemma
4.2. Below we assume that $\partial B\subset \partial H$.
Furthermore we assume that

(*) $|B\cap(F_{1}\cup F_{2})|$ is minimal among all compressing
disks $B$ of $\partial H_k$.

Suppose first $B\subset M_{k}$. Since $B\cap (E_1\cup
E_2)=\emptyset$,   $B$ separates a 3-ball $B^3$ from $J_k$ by
Lemma 2.3 (4). Since $\partial B$ is non-trivial in $\partial
H_K$, then either $ B^3$ contains only one of $ E_{1}$ and $
E_{2}$, see Figure 4.1 (c), and then $||a_i, B||\neq \emptyset$
for all $a_i\subset J_k$, a contradiction; or $B^3$ contains both
$E_{1}$ and $E_{2}$ and $k=2$ in this case, see Figure 4.1 (d),
$a_4$ and $a_8$ are properly homotopic in $(B^3, E_1\cup
E_2)\subset (J_k, E_1\cup E_2)$, which contradicts Lemma 2.3 (2).

Suppose then $B\cap (F_{1}\cup  F_{2})\neq \emptyset$. By Lemma
4.1 and the assumption (*), $B\cap(F_{1}\cup F_{2})$ consists of
arcs. Then an outmost arc $a$ of $B\cap (F_{1}\cup F_{2})\subset
B$ separates a disk $B_{0}$ from $B$ with $B_0\subset M_k$ for
some $k$. By (*) $a$ must be non-trivial in $F_{1}\cup F_2 $,
(otherwise $| B\cap (F_{1}\cup F_{2})|$ can be reduced by pushing
$B_0$ to a suitable side). Since $|B_0\cap (E_1\cup E_2)|=1$,
$B_{0}$ separates a 3-ball from $J_{k}$ by Lemma 2.3 (4). Then we
reach a contradiction by the same reason in the end of the proof
of Lemma 4.1 \qed

\vskip 0.3 true cm

{\bf Lemma 4.4.} {\it $M$ is anannular.}

{\bf Proof.} \ Suppose $H_{K}$ contains an essential annulus $A$.
Assume that

(**) $| A\cap (F_{1}\cup F_{2})|$ is minimal among all essential
annuli in $H_{K}$.

By Lemma 4.1 and (**), each component of $A\cap (F_{1}\cup F_{2})$
is non-trivial in both $A$ and $(F_{1}\cup F_{2})$. There are
three cases:

Case 1.  $\partial A\subset T$.  There  are two sub-cases:

(1)  $||\partial_1 A,\mu||=0$. May assume $\partial A\cap
(F_{1}\cup F_{2})=\emptyset$.

Suppose first $ A\cap (F_{1}\cup F_{2})=\emptyset$. May assume
that $A$ is contained in $M_{2}$ (the remaining  cases are more
directly). The whole $\partial A$ must lie in the same $A_i\subset
M_2$ by Lemma 2.3 (3). Since $A$ is essential, $A\cap
D_{i}=\emptyset$ for $i=2,4,8$. Then $\partial A\subset A_6$. Let
$M'$ be obtained by cutting $M_{2}$ along $ D_2, D_4, D_8$. Then
$A$ is still an incompressible annulus in $M^{'}$, and $D_6\subset
M_2$ become a properly embedding disk $D'_6\subset M^{'}$ with
$\partial D'_{6}\cap A_6=\partial D_{6}\cap A_6$, a non-trivial
arc of $A_6$. Since $\partial A\subset A_6$, $A\cap D'_{6}$ is an
arc in both $A$ and $D_{6}^{'}$. Hence there is a
$\partial$-compressing disk of $A$ in $M'$ which is also a
$\partial$-compressing disk of $A$ in $M_{2}$, which contradicts
that $A$ is  essential in $M_2$.

Suppose then $A\cap (F_{1}\cup F_{2})\neq \emptyset$, which  must
be a union of circles. An outmost circle $a$ of $A\cap (F_{1}\cup
F_{2})\subset A$ and $\partial_1 A$ bound an annulus $A^*\subset
A$. May assume that $a\subset F_{1}$, and $\partial_{1} A \subset
A_{i}\subset M_k$. Let $B^{*}$ be the disk bounded by $a$ on
$E_{1}$ and $D$ be the meridian disk of $N(K)$ bounded by
$\partial_1 A$. $B^{*}\cup A^*\cup D$ is a separating 2-sphere
$S^{2}$ which bounds a 3-ball $B^3\subset J_k$, see Figure 4.1
(e). Since $|(A^*\cup D)\cap (\cup_{a_j\subset J_k}
a_j)|=|(A^*\cup D)\cap  a_i)|=1$, by applying Lemma 2.3 (2) as
before we have $|B^*\cap (\cup_{a_j\subset J_k} a_j)|=|B^*\cap
a_i|=1$. Hence $a$ and $\partial_1 A_i$ bound an annulus $A'$ in

$F_1$. Now by pushing the annulus $\overline {A-A^*}\cup A'$
 to a suitable side of $F_1$, $|A\cap (F_{1}\cup
F_{2})|$ is reduced, which contradicts (**).

(2)  $||\partial_1 A, \mu||\geq 1$. Now $A\cap (F_{1}\cup F_{2})$
consists of non-trivial arcs in $A$,   which  cut $A$ into
$8||\partial_1 A, \mu||$ rectangles and each rectangles has two
opposite edges on $F_1\cup F_2$ and two opposite edges on $A_i$
and $A_{\pi(i)}$, where $\pi(i)=i+l \text{\, mod \,8}$. If $l=0$,
then the two ends of each arc of $A\cap (F_{1}\cup F_{2})$ lie in
a same component of $\partial(F_1\cup F_2)$, and an inner most arc
is trivial
 in $F_1\cup F_2$,
 a contradiction [G]. If $l\neq 0 \text{\,mod\,} 8$.
 Then $a_{6}$ and
$a_{6+l}$ are properly isotopy in $M_2$, which contradicts to
Lemma 2.3 (2).

Case 2. $\partial_{1} A\subset T$ and $\partial_{2}A\subset
\partial H$.

By Lemma 4.3 both $\partial H$ and $T$ are incompressible  in
$H_K$. Clearly $H_K$ is not homeomorphic to $T\times I$. Since
both Dehn fillings along $\mu$ and $\partial A_1$ compress
$\partial H_K$, by [2.4.3 CGLS], $\Delta (\partial_{1} A,\mu)\leq
1$. There are two sub-cases.

(1)  $||\partial_1 A,\mu||=0$. Since $\partial_{1} A$ is disjoint
from $F_{1}\cup F_{2}$, $\partial_{2} A$ is disjoint from
$F_{1}\cup F_{2}$ (otherwise there is an arc in $A\cap (F_{1}\cup
F_{2})$ with two ends in $\partial _2 A$ which is trivial in $A$).
Then it follows $A$ is disjoint from $F_{1}\cup F_{2}$ by the
proof of Case 1. Suppose $\partial_1 A\subset A_i\subset M_k$ for
some $i, k$. Let $D$ be the meridian disk of $N(K)$ bounded by
$\partial_1 A$ and $B=A\cup_{\partial_{1} A} D$. Then $B$ is a
proper disc in $J_k$, $\partial B$ is non-trivial in $\partial
H\cap J_k$, and $|B, \cup _{a_j\subset M_k} a_j|=|B, a_i|=1$.
Since $B\cap (E_1\cup E_2)=\emptyset$, $B$ separates a 3-ball
$B^3$ from $J_2$ by Lemma 2.3 (4). If $ B^3$ contains only one of
$E_{1}$ and $E_{2}$, then $B$ meets all $a_j$ in $J_k$. If $B^3$
contains both $E_{1}$ and $E_{2}$, then $B$ meets $a_i$ in
non-zero even number. In each case we reach a contradiction.

(2) $||\partial_{1} A, \mu||=1$. It is easy to see that this  case
is ruled out by Lemma 2.3 (5).

Case 3. $\partial A\subset \partial H$.

Suppose first $A\cap (F_{1}\cup F_{2})=\emptyset$. Since $A$ is
essential, $A$ is disjoint from $D_{i}$ for $i\neq 6$. May assume
that $A\subset M_{2}$ (the remaining  cases are the same). Since
$\partial A\subset \partial H\cap  J_2$ and $A$ is disjoint from
$c_{4}, c_{2}$, $A$ separates $J_2$ into two solid torus $J^*$ and
$J'$ such that $(E_1\cup E_2)\cap J'=\emptyset$ and separates a
disc $B'\subset J'$ from $B_2$, see Figure 4.1 (f). Since $J'$ is
disjoint from all $A_i\subset M_2$, $A$ is $\partial$-compressible
in $M_2$, which contradicts that $A$ is essential in $M_{2}$.

Suppose then $A\cap (F_{1}\cup F_{2})\neq \emptyset$. There are
two sub-cases:

(1) $A\cap (F_{1}\cup F_{2})$ consists of circles. Then an outmost
circle $a$ of $A\cap (F_{1}\cup F_{2})\subset A$ and $\partial_1
A$ bound an annulus $A^*\subset A$ such that $A^*\subset M_k$. Let
$B^{*}$ be the disk bounded by $a$ on $E_{1}\cup E_2$ and
$D^*=A^*\cup B^*$. By a slightly pushing, we have $D^*\subset
J_k$,  moreover (i) $D^*\cap a_i\ne \emptyset$ for some
$a_i\subset J_k$, (ii) for each $a_j\subset J_k$, $|D^*\cap
a_j|\le 2$, and $\le 1$ if $k=2$. Since $\partial D^* \subset
J_k\cap\partial H$, $D^*$  separates a 3-ball $B^3$ from $J_k$ by
Lemma 2.3 (4). Now either $B^3$ contains both $E_1$ and $E_2$,
$k=2$ in this case, and $D^*$ meets each $a_j\subset J_2$ in even
number of points, which contradicts (i) and (ii) above; or $B^3$
contains only one $E_i$, say $E_1$, then $\partial_1 A=\partial
D^*$ is parallel to $\partial E_1$. Since $A^*$ is disjoint from
$K$, $|B^*\cap (\cup_{a_j\subset J_k} a_j)|=|D^*\cap
(\cup_{a_j\subset J_k} a_j)|=|E_1 \cap (\cup_{a_j\subset J_k}
a_j)|=4$. Hence  $a$ and $\partial E_1$ bound an annulus $A'$ in
$F_1$ by applying Jordan Curve Theorem. Now we reach a
contradiction by the same argument the end of Case 1 (1).

(2) $A\cap (F_{1}\cup F_{2})$ consists of arcs. Then $F_{1}\cup
F_{2}$ cut $A$ into  rectangles $R_{i}$, and  each $R_i$ has two
opposite sides in $F_1\cup F_2$ and remaining two sides in
$\partial H$. Then $ R_{i}$ separates a 3-ball $B^{3}_{i}$ from
$J_k$ by Lemma 2.3 (4). Let $D_{i}^{1}$ and $D_{i}^{2}$ be two
disks of $B^3_i\cap (E_1\cup E_2)$. By Lemma 2.3 (2) we have (i)
$\partial_1 a_{j}\subset D_{i}^{1}$ if and only if $\partial_2
a_{j}\subset D_{i}^{2}$. (ii) $a_{j}$ and $a_{l}$ are contained in
the same $B^{3}_{i}$ implies that $(j,l)=(2,8)$.

If $a_2$ and $a_8$  belong to the same $B^3_i$, then so do  $a_3$
and $a_7$, which contradicts (ii). Hence there is only one $a_{j}$
in each $B^{3}_{i}$. Hence $A$ separates from $H$ a solid torus
with $K$ as centerline (up to isotopy). Then $K$ and  a component
of $\partial A\subset
\partial H$ bound an annulus, which contradicts Lemma 2.3 (5).\qed

\vskip 0.5 true cm

{\bf  \S 5. $H_{K}$ contains no closed essential surfaces.}

Recall $W$, $W_i$, $l_i$ defined in \S 2.

Suppose $H_{K}$ contains closed essential surfaces $F$. We define
the complexity of  $F$  by an ordered pair
$$C(F)=(|F\cap W|, |F\cap (F_{1}\cup F_{2})|).$$
Suppose $F$ realizes the minimality of $C(F)$. By the minimality
of $C(F)$, Lemma 4.1 and the standard argument in 3-manifold
topology, we have

\vskip 0.3 true cm

{\bf Lemma 5.0.}

{\it (1) each component of $F\cap(F_{1}\cup F_{2})$ is a
non-trivial circle in both $F$ and $F_{1}\cup F_{2}$,

(2) $F\cap W \subset W$ is a union of arcs as in Figure 5.1. Hence
$|F\cap l_{i}|=|F\cap l_{j}|$ for all $i,j$.}

\begin{center}%
\includegraphics[totalheight=4cm]{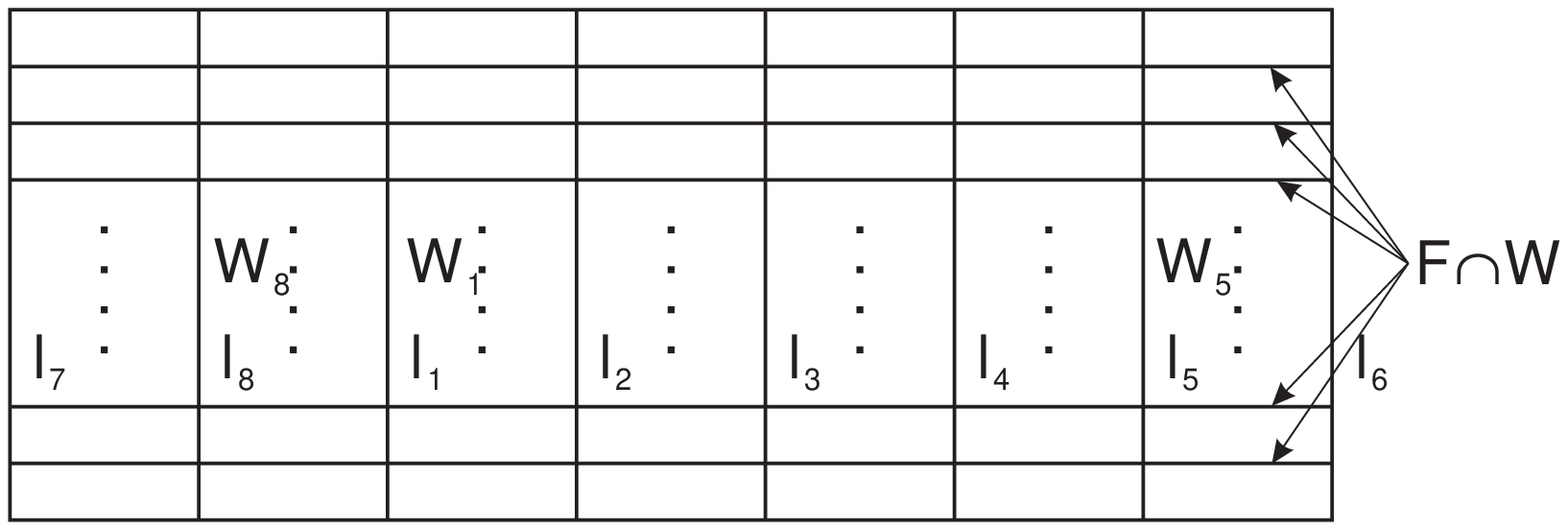}%
\begin{center}%
Figure 5.1
\end{center}
\end{center}

The positions of $\partial A_{i}$ and $l_i$ in $F_1\cup F_2$ are
indicated as in Figure 5.2.

\vskip 0.3 true cm

\begin{center}%
\includegraphics[totalheight=4.5cm]{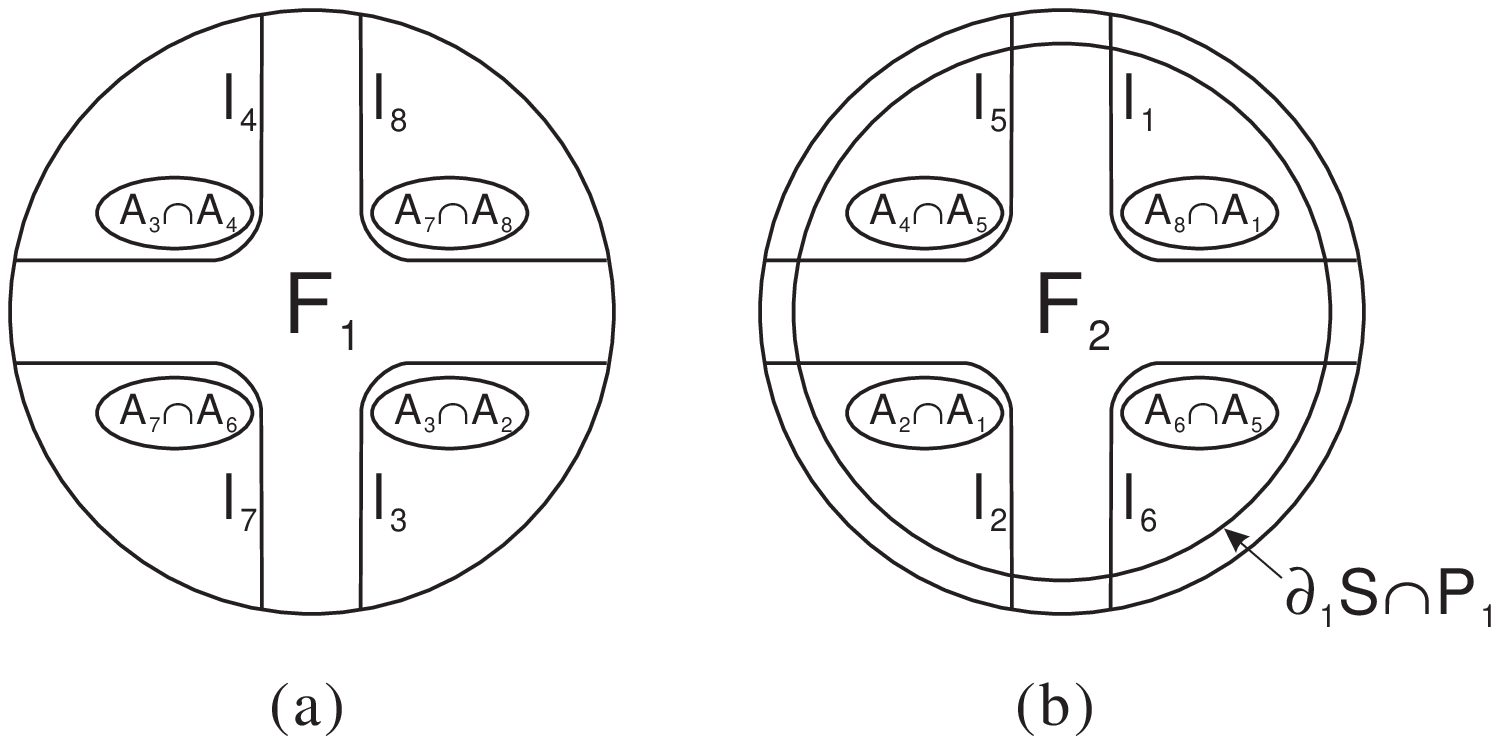}%
\begin{center}%
Figure 5.2
\end{center}
\end{center}

Below we will use $\tilde s$ to denote a given family of parallel
disjoint proper 1-manifolds on some surface, and use $s$ to denote
a representive (a component) of $\tilde s$.

{\bf Lemma 5.1} {\it Each component of $F\cap M_{k}$  is isotopic
to either $ M_{k}\cap\partial H$ or some $A_i\subset M_k$, where
$k=1,3$.}

{\it Proof.} The proofs  for $k=1$ and $3$ are the same. Assume
$k=3$. First we need

{\bf Lemma 5.2} {\it Components of $F\cap  F_2$ in $F_2$ which are
not parallel to a component of $\partial F_2$ are divided into two
families of circles  $\tilde s_1$ and $\tilde s_2$ in Figure 5.3
(a) and (c). Moreover in each case $|\tilde s_1|=|\tilde s_2|$.}

\vskip 0.3 true cm

\begin{center}%
\includegraphics[totalheight=4.5cm]{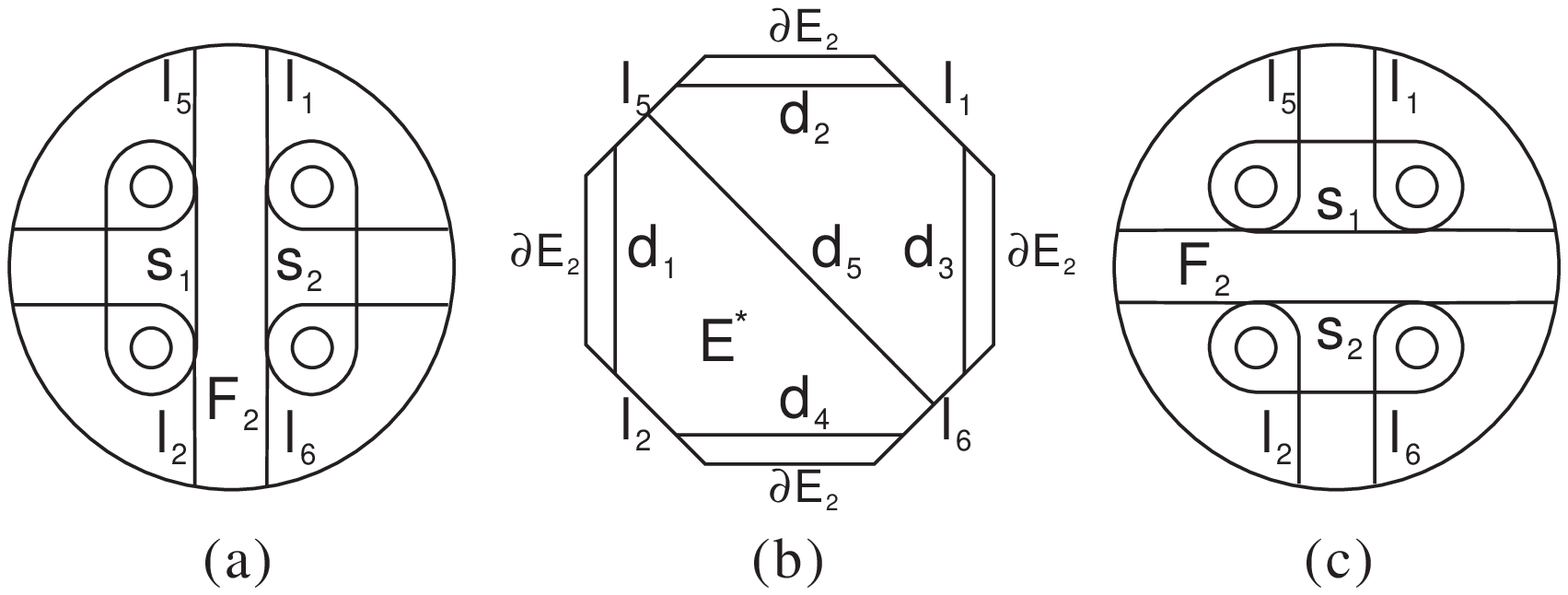}%
\begin{center}%
Figure 5.3
\end{center}
\end{center}

{\bf Proof.}  Since each component of $F\cap F_2$ isotopic to a
component of $\partial F_2$ contributes the same to $|F\cap l_k|$
 for all $l_k\subset F_2$, we may assume that $F\cap F_2$
contains no such components when we apply Lemma 5.0 (2) to prove
Lemma 5.2.

Note that $l_1, l_2, l_5, l_6$ separate $F_2$ into four annulus
and one disc which is presented as an octagon  $E^*$ in Figure 5.3
(b), where $F\cap E^*$ are presented as five families of proper
disjoint arcs $\tilde d_1, ..., \tilde d_4, \tilde d_5$ with
$\partial d_i$ in different $l_j$ and $l_k$ for each $i$. By Lemma
5.0 (2), we have

$|\tilde d_4|+|\tilde d_1|=|\tilde d_2|+|\tilde d_1|+|\tilde
d_5|=|\tilde d_2|+|\tilde d_3|=|\tilde d_4|+|\tilde d_3|+|\tilde
d_5|.$

It follows $|\tilde d_5|=0$, $|\tilde d_1|=|\tilde d_3|$ and
$|\tilde d_2|=|\tilde d_4|$. Back to Figure 5.2 (b), since no
component is isotopic to $\partial F_2$, it follows that either
$|\tilde d_2|=0$, which is Figure 5.3 (a), or $|\tilde d_1|=0$,
which is Figure 5.3 (c). \qed

Let us return to  the proof of Lemma 5.1. Let $S$ be a component
of $F\cap M_3$. Each $W_{i}$ separates a solid tori $P_i$ from
$M_3$, $i=1,5$. Let $M^{'}_3=\overline{M_3-(P_1\cup P_5)}$, which
is a solid torus. There are three cases to discuss.

Case 1.  If a component of $\partial S$ is isotopic to a component
of $\partial A_{i}$, $i=1,5$. By the minimality of the complexity
$C(F)$, $S$ is disjoint from $W_{i}$, and therefore $S\subset
P_i$, which is an annulus isotopy to $A_{i}$, $i=1$ or $5$.

Case 2.  If a component of $\partial S$ is isotopic to $\partial
E_{2}$, let $\partial_{1} S$ be the outmost component of $\partial
S\subset  F_2$ which is isotopic to $\partial E_{2}$. Now
$\partial _{1} S$ intersects $P_{i}$ as in Figure 5.2 (b) and
$W_i\cap S$ contains two arcs $b_i^*$ and $b_i^\#$ with ends in
$\partial _1 S$. Let
 $S_i$ be the component of $S\cap P_i$
which meets $\partial _1 S$, $i=1,5$.  Then $S_i$ is
incompressible in $P_i$. Since $\partial S_i=(\partial_1 S\cap P_i
)\cup (b_i^*\cup b_i^\#)$ bounds a disk in $P_{i}$ parallel to
$\partial M_{3}$, $S_i$ itself is such a disc, $i=1,5$.  Let $S_3$
be a component of $S\cap M'_3$ which meets $\partial_1 S$, then
$S_3$ is incompressible in the solid torus $M'_3$ and $\partial
S_3$ has a component $(\partial _1 S\cap M_3')\cup (b_1^*\cup
b_1^\#)\cup (b_5^*\cup b_5^\#)$ which bounds a disk in $\partial
M_{3}'$ as in Figure 5.4. Hence $S_3$ itself is such a disk. Thus
$S=S_1\cup_{b_1^*\cup b_1^\#} S_3\cup _{ b_5^*\cup b_5^\#} S_5$ is
isotopic to $M_{3}\cap\partial H$.

\vskip 0.3 true cm

\begin{center}%
\includegraphics[totalheight=4.5cm]{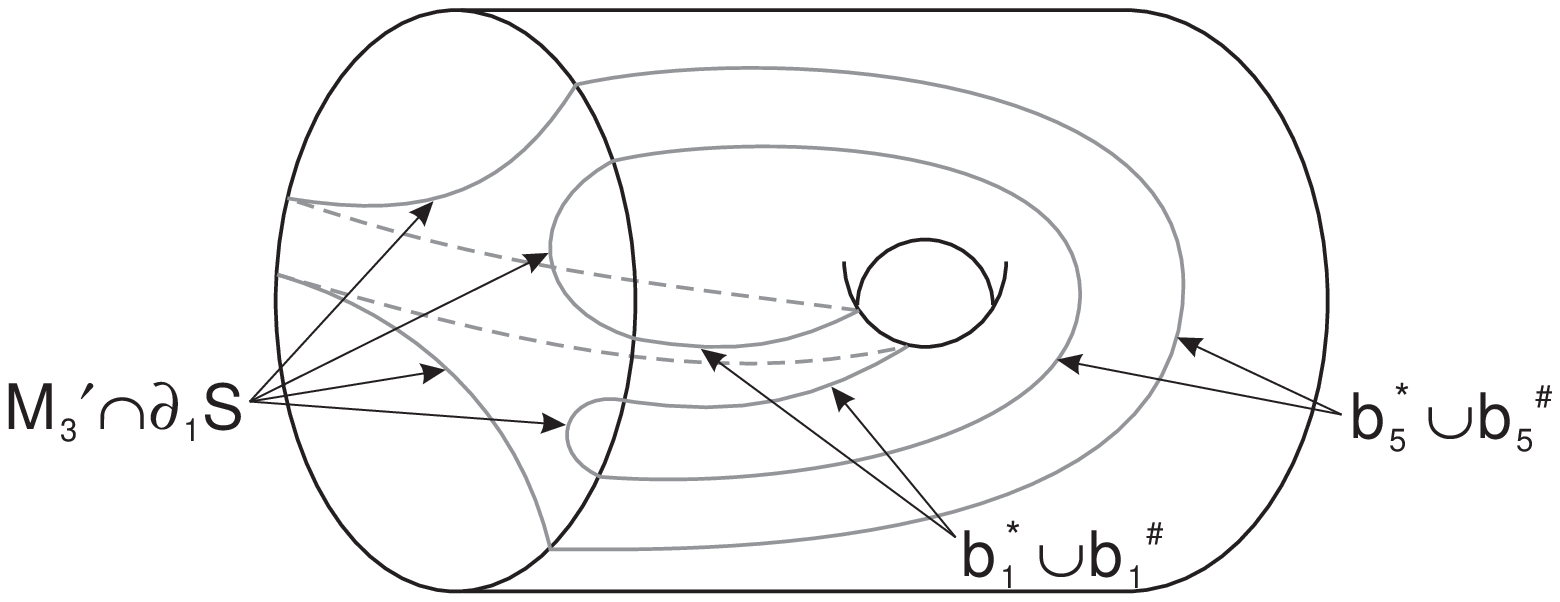}%
\begin{center}%
Figure 5.4
\end{center}
\end{center}

\vskip 0.3 true cm

\begin{center}%
\includegraphics[totalheight=4.5cm]{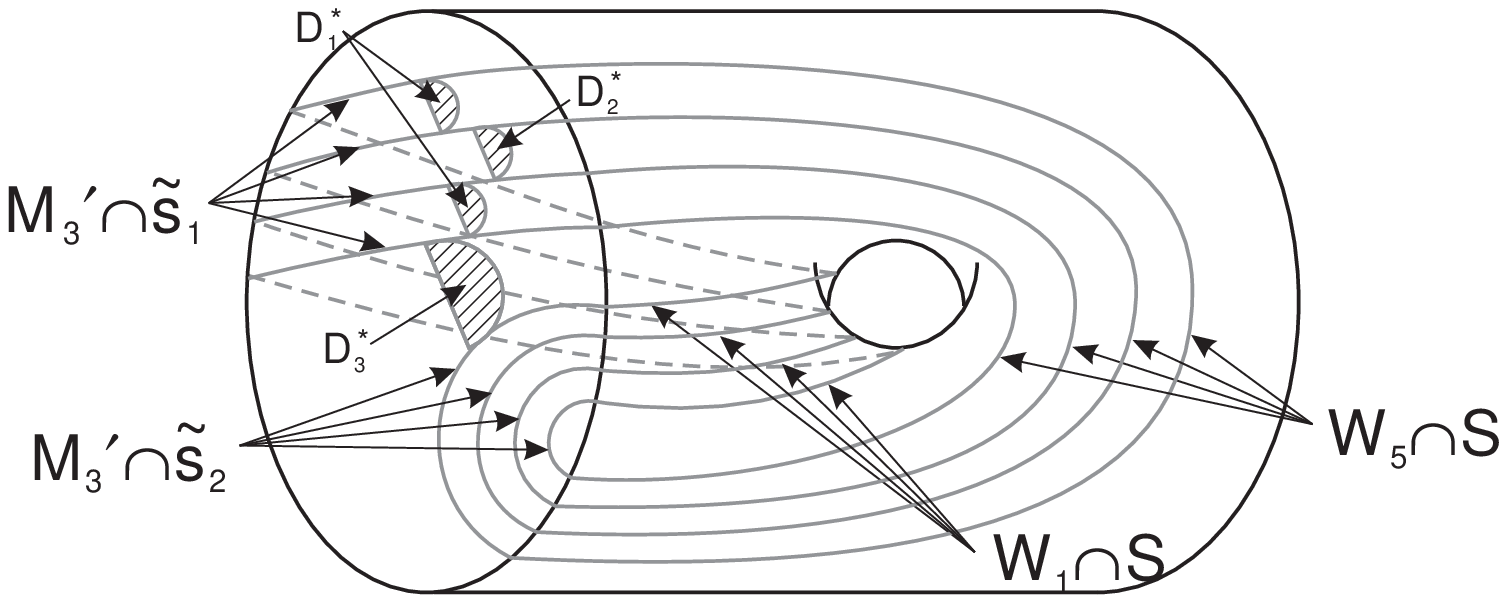}%
\begin{center}%
Figure 5.5
\end{center}
\end{center}

By Lemma 5.2, to finish the proof we need only to rule out  Case 3
below.

Case 3. $|\tilde s_{1}|=|\tilde s_{2}|\neq 0$ in Figure 5.3 (a) or
(b). Since the discussion for (a)  and (c) in Figure 5.3 are the
same, we just discuss the former case.

We may assume that no component of $S$ is isotopic to a component
of $\partial F_2$ by Case 1 and  Case 2 we just discussed. Let
$S_3'=S\cap M'_3$. Then $\partial S^{'}_3$ contains $2|\tilde
s_1|$ circles which are produced from the arcs $(\tilde s_1\cup
\tilde s_2) \cap M'_3$ and the arcs $(W_1\cup W_5)\cap S$, as in
Figure 5.5, where $|\tilde s_1|=2$. Note each circle in $\partial
S'_3$ is non-trivial in $\partial M'_3$. Since $S_3'$ is
incompressible in the solid torus $M_3'$, each component of
$S^{'}_3$ is an annulus which is $\partial$-compressible. Now
$B_3'\cap S_3'$ is a union of arcs, where $B_3'=B_3\cap M_3'$. An
outmost arc $b$ of $B_3'\cap S_3'\subset B_3'$ separates a disc
$D^*$ from  $B_3'$. As a $\partial$-compressing disc of $S'_3$,
$D^*$ can be moved into the positions of $D^*_1, D^*_2, D^*_3$,
indicated as in Figure 5.5. Now back to $M_3$, those $D^*_i$'s in
Figure 5.5 are corresponding to those $D^*_i$'s in Figure 5.6 (a),
$i=1,2,3$. In the cases of $D^*_1$ and $D^*_3$ in Figure 5.6 (a),
one can push $F$ along the disc to reduce $|F\cap W|$; in the case
of $D^*_2$ in Figure 5.6 (a), one can push $F$ along the disc to
reduce $|F\cap (F_{1}\cup F_{2})|$, but not to increase $|F\cap
W|$. In each case, it contradicts the minimality of $C(F)$. \qed

{\bf Remark on Figures 5.4, 5.5 and 5.7.} In Figure 5.4, to
simplify the picture, $W_1$ does not meet $B_3$ in three arcs as
it should be. But one verifies easily that this simplification
does not affect the proof. The same remark is needed for Figures
5.5 and 5.7. Moreover in Figure 5.7, we only drawn a representive
$e_i$ for a families $\tilde e_i$ and so on.

\vskip 0.3 true cm

\begin{center}%
\includegraphics[totalheight=4.5cm]{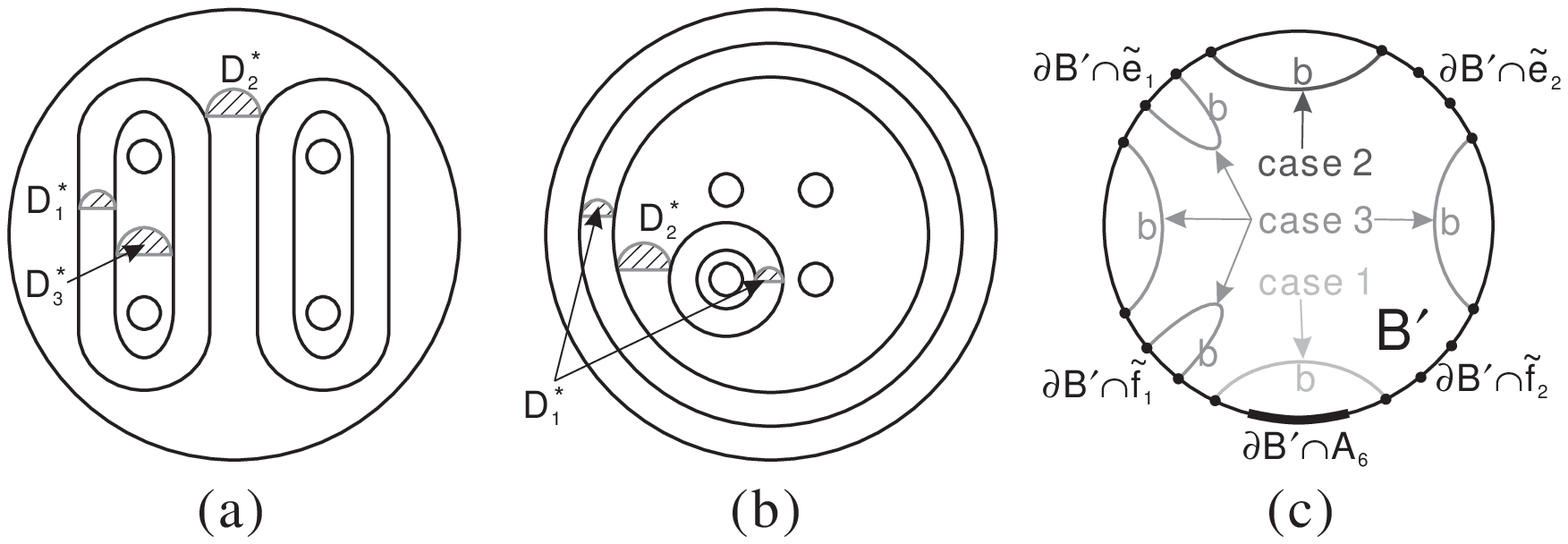}%
\begin{center}%
Figure 5.6
\end{center}
\end{center}

{\bf Proposition 5.3.} {\it $H_{K}$ contains no closed essential
surface.}

{\bf Proof.}  Now we consider $F\cap M_2$.  Each component of
$F\cap (F_1\cup F_{2})$ is isotopic to a component of $\partial
F_{2}\cup
\partial F_{1}$ by Lemma 5.1. Apply Lemma 5.0 (2) again, we have $|\widetilde{\partial
E_2}|=|\widetilde  {\partial E_1}|$, where $\widetilde  {\partial
E_i}\subset F_i$ are components of $F\cap F_i$ isotopic to
$\partial E_i$, $i=1,2$. Each $W_{i}$ separates a solid tori $P_i$
from $M_2$, $i=2,4,8$. Let $M^{'}_2=\overline{M_2-(P_2\cup P_4\cup
P_{8})}$, which is a handlebody of genus 2.

Note if $F\cap M_2$ has a component $S$ such that a component of
$\partial S$ is isotopic to $\partial A_i$, $i=2,4,8$, by the
minimality of $C(F)$, $S\subset P_i$ and hence $S$ is isotopic to
$A_i$, $i=2,4,8$.

Let $S^{'}_2=F\cap M_2^{'}$. Then  $\partial S^{'}_2$  consists of
possibly five families of circles $\tilde e_1, \tilde e_2, \tilde
e_{3}$, $\tilde f_1$ and $\tilde f_2$, where $\tilde e_1$, $\tilde
e_2$ and $\tilde e_3$ with $|\tilde e_i|=|\widetilde {\partial
E_1}|$ are produced from  the arcs $(\widetilde {\partial E_1}\cup
\widetilde {\partial E_2}) \cap M'_2$ and the arcs of $(W_2\cup
W_4\cup W_8)\cap F$ with end points lying $\widetilde {\partial
E_1}\cup\widetilde  {\partial E_2}$, $\tilde f_1\subset F_1$ and
$\tilde f_2\subset F_2$ are parallel copies of the two components
of $\partial A_6$ respectively. All those are indicated in Figure
5.7 (see Remark on Figures 5.4., 5.5, 5.7). Moreover

(i) each component of $\tilde e_{3}$ bounds a disk in $\partial
M_2'$, hence bounds also a disk in $F$.

(ii) any two components in $\tilde e_1\cup \tilde e_2$ bound an
annulus in $\partial M_2'$ disjoint from $A_6$.

\vskip 0.3 true cm

\begin{center}%
\includegraphics[totalheight=5cm]{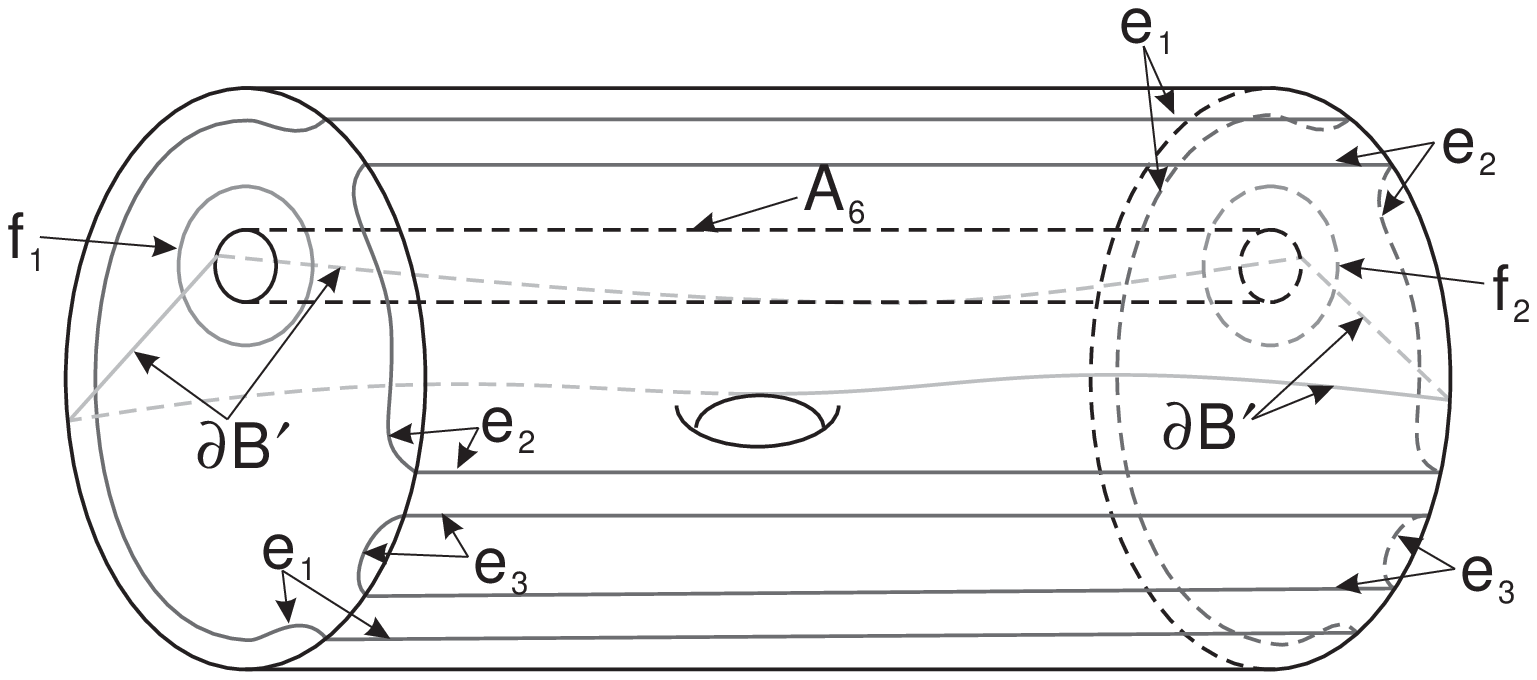}%
\begin{center}%
Figure 5.7
\end{center}
\end{center}

There is a proper disc $B^{'}$ in $M_2'$ with $\partial B'$ shown
 in Figure 5.7 such that

(iii) $\partial B'$ meets those four families in the cyclic order
$\tilde e_{1}$, $\tilde e_{2}$, $\tilde f_2$, $\tilde f_1$,

(iv) $\partial B'$ meets each component of  $\tilde e_{1}\cup
\tilde e_{2}\cup \tilde f_1\cup \tilde f_2$ in one point and
$\partial B'\cap A_6$ is a non-trivial arc in $A_6$,

Let $S'$ be a component of $S'_2$. Since $S^{'}$ is incompressible
in $M_2'$,  $S^{'}\cap B^{'}$ consists of arcs. By (iv) there is
an outmost arc $b$ of $S' \cap B'\subset B'$  which separates a
 disk $D^*$ from  $B'$ so that

(v) $\partial D^*$ disjoint from $A_6$ and $D^*$ is a
$\partial$-compressing disk of $S'$.

We divide the remaining discussion into three cases by (iii).
(Figure 5.6 (c) is helpful to understand (iii), (iv) and (v) above
and each case below.)

Case 1.  One end of $b$ is in $f_1\in \tilde f_1$ and the other is
in $f_2\in \tilde f_2$. In this case $|\widetilde {\partial
E_1}|=0$ and all $\tilde e_i$'s do not exist by (iii) and (v). One
can show that $S'$ is isotopic to $A_6$ by cutting and pasting
argument in 3-manifold topology, the detail is contained in what
we will do in Case 2.

Case 2.   One end of $b$ is in $e_1\in \tilde e_1$ and the other
is in $e_2\in \tilde e_2$. By (ii), $e_1$ and $e_2$ bound an
annulus $A$ in $\partial M_2'$ disjoint from $A_6$.  Let $M_2''$,
$S''$,  $A'$, $A_6'$,  $e_1'$, $e_2'$  be the images of $M_2'$,
$S'$, $A$, $A_6$, $e_1$, $e_2$ respectively after cutting $M_2'$
along $B'$. By (iv), $e_i'$ is an arc, $i=1,2$. Let $b_i$,
$D^*_i$, $i=1,2$, be the two copies of $b$ and $D^*$ after cutting
$M_2'$ along $B'$. By (v), the circle $c'\subset \partial S''$
formed by four arcs $e'_1$, $e_2'$, $b_1'$ and $b_2'$  bound a
disc $D^*_1\cup A'\cup D^*_2$ in $\partial M_2''$ which is
disjoint from $A_6'$. Since $S''$ is incompressible in $M_2''$,
$S''$ is such a disc up to isotopy. Back to $M_2'$, $S'$ is
isotopic to the annulus $A\subset M_2'$. Back to $M_2$, by (i) and
similar argument in Case 2 in the proof of Lemma 5.1, $S$ is
isotopic to $M_{2}\cap\partial H$.

Case 3.  If either $\partial b$ lie in one of the four families
$\tilde e_1$, $\tilde e_2$, $\tilde f_1$ and $\tilde f_2$, or  one
end of $b$ is in $\tilde e_i$ and the other  in $\tilde f_i$,
$i=1$ or 2,  then $D^*$ can be moved in $M_2'$ keeping to be a
$\partial$-compressing disk of $S'$ so that when we go back to
$M_2$ it is a $\partial$-compressing disk of $F\cap M_2$ in the
position of either $D^*_1$ or $D^*_2$ in Figure 5.6 (b).  One can
push $F$ along either $D_1^*$ or $D^*_2$ to reduce $C(F)$, which
contradicts the minimality of $C(F)$. (Refer the end of the
argument in Case 3 of the  Proof of Lemma 5.1).

 So  each component $S$ of $F\cap M_{2}$ is isotopic to either $
M_{2}\cap\partial H$ or $A_i$, $i=2,4,6,8$. In the former case,
$\partial S$ is $\partial E_1$ and $\partial E_2$ which bound (up
to isotopy) $\partial H\cap M_1$ and $\partial H\cap M_3$
respectively by Lemma 5.1, and then $F$ is isotopic to $\partial
H$. In the later case,  by Lemma 5.1, each component of $F\cap
(M_{1}\cup M_{3})$ is an annulus isotopic to one of $A_{1}, A_{3},
A_{5}, A_{7}$. Since $F$ is closed,  it follows that $F$ is a
torus isotopic to $T$.\qed

Proposition 3.0 follows from Lemmas 4.2, 4.3, 4.4 and Proposition
5.3. Hence Theorem 1 is proved.

{\bf References.}

\vskip 0.4 true cm

[CGLS] M.
 Culler, C. Gordon, J. Luecke and P. Shalen, Cyclic
surgery on knots, Ann. of Math., 125(1987), 237-300.

[G] C. Gordon, Dehn filling: a survey. Knot theory (Warsaw, 1995),
129--144, Banach Center Publ., 42, Polish Acad. Sci., Warsaw,
1998.

[H] A. Hatcher, On the boundary curves of incompressible surfaces.
Pacific J. Math. 99 (1982), 373-377.

[J] W. Jaco, Adding a 2-handle to  a 3-manifold, An application to
Property R, Proc. AMS 92 (1984), 288-292.

[L] M. Lackenby, Attaching handlebody to 3-manifolds, Geometry and
Topology, Vol. 6 (2002), 889-904 (2002)

[Q] R.F. Qiu,  Incompressible surfaces in handlebodies
    and closed 3-manifolds of
      Heegaard genus two, Proc. AMS 128(2000), 3091-3097.

[SW] M. Scharlemann and Y. Wu, Hyperbolic manifolds and
     degenerating handle additions. J. Aust. Math. Soc. (Series A) 55
     (1993), 72-89.

[T] W. Thurston, Three dimensional manifolds, Kleinian groups and
    hyperbolic geometry. Bull. AMS, Vol. 6, (1982) 357-388.

\vskip 0.4 true cm

Ruifeng Qiu

Department of Mathematics,

Dalian University of Technology, 130023, China

qiurf@dlut.edu.cn

\vskip 0.4 true cm

Shicheng Wang

 Department of Mathematics,

Peking University, 100871, China

 wangsc@math.pku.edu.cn

\enddocument
\bye